\pdfoutput=1
\documentclass[a4paper,12pt]{article}
\usepackage{amsmath, amsthm, amssymb}
\usepackage{epsfig}
\usepackage{color}
\usepackage{comment}
\usepackage{mathdots}
\usepackage{mathrsfs}
\usepackage{epstopdf}
\usepackage{graphicx}
\usepackage{makeidx} 
\usepackage{ifpdf}
\usepackage[bookmarks=false]{hyperref}

\def\IR{{\mathbb R}}
\def\IC{{\mathbb C}}

\def\IL{{\mathbb L}}

\newcommand{\sIL}{{{{\mathbb L}_s}}}

\newtheorem{remark}{Remark}[section]

\newcommand{\bA}{{\bf A}}
\newcommand{\bB}{{\bf B}}
\newcommand{\bC}{{\bf C}}
\newcommand{\bD}{{\bf D}}
\newcommand{\bE}{{\bf E}}

\newcommand{\bS}{{\bf S}}
\newcommand{\bY}{{\bf Y}}

\newcommand{\bL}{{\bf L}}
\newcommand{\bM}{{\bf M}}

\newcommand{\bI}{{\bf I}}
\newcommand{\bH}{{\bf H}}

\newcommand{\bW}{{\bf W}}
\newcommand{\bR}{{\bf R}}

\newcommand{\bX}{{\bf X}}
\newcommand{\bx}{{\bf x}}
\newcommand{\by}{{\bf y}}
\newcommand{\bu}{{\bf u}}
\newcommand{\bV}{{\bf V}}
\newcommand{\bU}{{\bf U}}

\newcommand{\bfe}{{\bf e}}

\newcommand{\bz}{{\bf z}}

\newcommand{\bfz}{{\mathbf 0}}
\newcommand{\cC}{ {\cal C} }

\newcommand{\cE}{ {\cal E} }
\newcommand{\cO}{ {\cal O} }

\newcommand{\cP}{ {\cal P} }

\newcommand{\balpha}{ \boldsymbol{\alpha} }

\newcommand{\bLambda}{\boldsymbol{\Lambda}}

\newcommand{\cR}{ {\cal R} }
\newcommand{\cT}{ {\cal T} }

\newcommand{\bLa}{\boldsymbol{\Lambda}}

\newtheorem{theorem}{Theorem}
\newtheorem{lemma}[theorem]{Lemma}

\begin{document}

\title{Rational approximation of the absolute value function from measurements: a numerical study of recent methods}

\author{I.V. Gosea\thanks{I. V. Gosea the is with Max Planck Institute for Dynamics of Complex Technical Systems, Magdeburg, Germany, (e-mail: gosea@mpi-magdeburg.mpg.de)} \ and A.C. Antoulas\thanks{A.C. Antoulas is with the Department of Electrical and Computer Engineering, Rice University, Houston, Baylor College of Medicine, Houston, and Max Planck Institute, Magdeburg, Germany, (e-mail: aca@rice.edu)}}

\date{\today}

\maketitle

\begin{abstract}
In this work, we propose an extensive numerical study on approximating the absolute value function. The methods presented in this paper compute approximants in the form of rational functions and have been proposed relatively recently, e.g., the Loewner framework and the AAA algorithm. Moreover, these methods are based on data, i.e., measurements of the original function (hence data-driven nature). We compare numerical results for these two methods with those of a recently-proposed best rational approximation method (the minimax algorithm) and with various classical bounds from the previous century. Finally, we propose an iterative extension of the Loewner framework that can be used to increase the approximation quality of the rational approximant.
\end{abstract}

\section{Introduction}\label{sec:intro}

Approximation theory is a well-established field of mathematics that splits in a wide variety of subfields and has a multitude of applications in many areas of applied sciences. Some examples of the latter include computational fluid dynamics, solution and control of partial differential equations, data compression, electronic structure calculation, systems and control theory (model order reduction reduction of dynamical systems).

 In this work we are concerned with rational approximation, i.e, approximation of given functions by means of rational functions. More precisely, the original to be approximated function is the absolute value function (known also as the modulus function). Moreover, the methods under consideration do not necessarily require access to the exact closed-form of the original function, but only to evaluations of it on a particular domain. Hence, these methods are data-driven.  

Since rational functions can be written as a ratio of two polynomials (in the numerator and in the denominator), they prove to be more versatile than simple polynomials. Consequently, rational functions can approximate functions with singularities and with oscillatory behavior better than polynomials can. Furthermore, rational functions are extensively used in systems and control theory since the input-output behavior of linear dynamical systems is described in the frequency domain by the so-called transfer function (which is a rational function). Finally, the location of the poles of transfer functions determines important system properties such as asymptotic stability, transient behavior or damping.

For more details of the many aspects of approximation theory, we refer the reader to the classical books in \cite{Da75,Po81}, and to the more recent books \cite{Tr13,Is18} that put also emphasis on practical issues such as software implementation of the methods and applications, i.e., on data analysis for the latter.

In this work we study two recent data-driven rational approximations methods: the Loewner framework introdued in \cite{ajm07}, and the AAA algorithm introduced in \cite{nst18}. The main contribution of this paper is to present an extensive numerical study of the performance of these two methods for approximating the absolute value function. We put an emphasis on the Loewner framework and show that this method yields, in some cases, very promising results (e.g., maximum approximation errors lower than $10^{-10}$). Moreover, we compare the numerical results with those of best rational approximants (computed by means of a recent numerical tool proposed in \cite{fntb18}). A multitude of bounds that were proposed in the last century are also included in the study. Finally, another main contribution of this work is to propose a novel Loewner iterative procedure that is able to further improve the approximation quality, and in some cases, to come near to the error provided by best rational approximants. 

This paper is organized as follows; after a brief introduction is given in Section\;\ref{sec:intro}, we continue with a historical account of approximating the absolute value function provided in Section\;\ref{sec:approx}. Afterwards, in Section\;\ref{sec:main}, we introduce the main rational approximation methods (the Loewner framework, the AAA and minimax algorithms) and we explain how the maximum approximation error is computed. Then, in Section\;\ref{sec:numerics}, we provide a thorough and extensive numerical study of the methods. In Section\;\ref{sec:iterate} the new Loewner iterative procedure is introduced (and also some numerical tests are provided), while in Section\;\ref{sec:conc} the conclusion follows.

\section{Approximation of the absolute value function - a brief historical account}\label{sec:approx}

The absolute value function $\vert x \vert$ is a continuous function, non-differentiable at $x = 0$. Approximating this function by polynomials played a significant role in the early development of approximation theory. In 1908, the Belgian mathematician de la Vall\'ee-Poussin raised the question of finding the best approximation of the absolute value function on the interval $[-1,1]$. This problem attracted the attention of several leading mathematicians of that period.

Approximating the function $\vert x \vert$ by means of polynomials or rational functions was studied intensively starting with the beginning of the 20th century. In what follows we present a short historical account of some of the main breakthrough results of quantifying the approximation error, either by means of polynomials or by rational functions. We follow the information provided in the survey papers \cite{Re87,vrc92} and in Chapter 4 of the book \cite{Fi03}. 

Denote the set of polynomials with real coefficients of degree at most n with $\cP_{n}$, where 
 \begin{equation}
\cP_{n} = \Big{\{} p(x) \Big{\vert} p(x) \in \mathbb{R}_n[x]  \Big{\}}.
\end{equation}
and with $\cE_n(\vert x \vert ; [-1,1])$ the error of the best polynomial approximation from $\cP_n$ of the function $\vert x \vert$ on the interval $[-1,1]$. More precisely, write
\begin{equation}\label{best_err_def}
\cE_n(\vert x \vert ; [-1,1]) := \inf\{ \Vert x - p(x) \Vert_{L_{\infty}[-1,1]} \  \vert \ p \in \cP_n \}.
\end{equation}
For brevity purposes, we sometimes use notation $\cE_n$ to refer to the quantity introduced in (\ref{best_err_def}).

One of the first significant results regarding the error for approximating the absolute value function by means of polynomials was presented in \cite{Be13}. There, it was shown that there exists a real positive constant $\beta$ (the so-called Bernstein constant) in the interval $[0.278,0.286]$ so that the following holds
\begin{equation}\label{Bernstein_ct}
\beta = \lim_{n \rightarrow \infty} 2n \cE_{2n}(\vert x \vert; [-1,1]).
\end{equation} 
Moreover, it was conjectured in \cite{vc85} that $2n \cE_{2n}$, which is a function of $n$, admits the following series expansion as $n \rightarrow \infty$
\begin{equation}
2n \cE_{2n}(\vert x \vert; [-1,1]) \approx \beta - \frac{K_1}{n^2} + \frac{K_2}{n^4} - \frac{K_3}{n^6}+\cdots, \ \ \ K_i >0. 
\end{equation} 
Other contributions were made throughout the years, such as the ones in \cite{Be52,Lo66}. The bottom line is that the error of the best polynomial approximation of the absolute value function is indeed $\cO(n^{-1})$ (as $n$ increases, the error decreases linearly with $n$).

Consequently, the result presented by D. J. Newman in \cite{Ne64} was indeed surprising. More precisely, he proved that the error of the best uniform approximation using rational functions of the function $\vert x \vert$  in the interval $[-1,1]$ is in $\cO(e^{-c \sqrt{n}})$. This meant that the speed of rational approximation is much faster than that of polynomial approximation.

Next, we introduce the following set of all real rational functions $r(x) = \frac{p(x)}{q(x)}$:
\begin{equation}
\cR_{m,n} = \Big{\{} \frac{p(x)}{q(x)} \Big{\vert} p(x) \in \mathbb{R}_m[x], \ \  q(x) \in \mathbb{R}_n[x]  \Big{\}}.
\end{equation}
Introduce the minimum approximation error in the following way:
\begin{equation}
\cE_{m,n}(\vert x \vert ; [-1,1]) := \inf\{ \Vert x - r(x) \Vert_{L_{\infty}[-1,1]} \  \vert \ r \in \cR_{m,n} \}.
\end{equation}
The result proven in \cite{Ne64} can be formulated as follows: for all $n \geqslant 4$, the following holds
\begin{equation}\label{Newman_bound}
\frac{1}{2}e^{-9 \sqrt{n}} \leqslant \cE_{n,n} (\vert x \vert; [-1,1]) \leqslant 3 e^{- \sqrt{n}}. 
\end{equation}
For more details on this result, and on some extensions of it to approximating other important functions, we refer the reader to \cite{Tu66,Go67} and to chapter 4 from the book \cite{PP87}.

Moreover, it was shown in \cite{Bu68} that the following inequality holds for all $n \geqslant 0$
\begin{equation}\label{Bulanov_bound}
 \cE_{n,n} (\vert x \vert; [-1,1]) \geqslant e^{- \pi \sqrt{n+1}},
\end{equation}
and in \cite{Vj75} that there exist positive constants $M_1$ and $M_2$ so that the following holds for all values of $n \geqslant 1$
\begin{equation}
M_1 \leqslant e^{\pi \sqrt{n}} \cE_{n,n} (\vert x \vert; [-1,1]) \leqslant M_2. 
\end{equation}
Finally, years later in the contributions \cite{St92-1,St92-2}, the following result was proven:
\begin{equation}\label{Stahl_bound}
\lim_{n \rightarrow \infty} e^{\pi \sqrt{n}} \cE_{n,n} (\vert x \vert; [-1,1]) =8. 
\end{equation}
The existence of the limit in (\ref{Stahl_bound})  (with exact value $8$) has been conjectured in \cite{vrc92} on the basis of extensive numerical calculations. The limit can be considered as the analogue of Bernstein's constant in (\ref{Bernstein_ct}) for the case of the best
polynomial approximation of $\vert x \vert$ on $[-1, 1]$.
To the best of our knowledge, the formula in (\ref{Stahl_bound}) represents the tightest possible attainable result.

\subsection{Newman's rational approximant}\label{sec:Newton}

In his groundbreaking paper from 1964 \cite{Ne64}, Newman not only provided bounds for the best rational approximation of $\vert x \vert$ but also showed how to compute a rational approximant explicitly. 

Let $n >0$ be a positive integer and let $\alpha = e^{-\frac{\sqrt{n}}{n}} \in (0,1)$ and $\epsilon = \alpha^n = e^{-\sqrt{n}}$. Introduce the polynomial $p(x)$ and the rational function $f(x)$ as
\begin{equation}
p(x) = (x+\alpha)(x+\alpha^2)\cdots(x+\alpha^{n-1}), \ \ \ \text{and} \ \ f(x) = \frac{p(-x)}{p(x)}.
\end{equation}
In \cite{Ne64} the following result was proven for all $x \in [\epsilon,1]$ and $n \geqslant 5$
\begin{equation}\label{ineq_f}
\vert f(x) \vert =  \Big{\vert} \frac{p(-x)}{p(x)} \Big{\vert}  \leq e^{\frac{2(\alpha-\alpha^n)}{\ln(\alpha)}}.
\end{equation}
To prove this, one uses the fact that $ \vert f(x) \vert  \leq  \vert f(1) \vert$ for $x \in [\epsilon,1]$ and that $\frac{1-\alpha^k}{1+\alpha^k} \leq e^{-2 \alpha^k} $. 

Now, choosing $n \geq 5$ it follows $2(\alpha-\alpha^n) >1$. Also note that $\ln(\alpha) = \ln(e^{-\frac{\sqrt{n}}{n}}) = -\frac{\sqrt{n}}{n}$.

 Hence, we have $\frac{2(\alpha-\alpha^n)}{\ln(\alpha)} < e^{-\sqrt{n}}$ and thus from (\ref{ineq_f}), it follows $\vert f(x) \vert < e^{-\sqrt{n}} \ \forall \ x \in [\epsilon,1]$ and $n \geqslant 5$. Finally, the order $(n-1,n-1)$ Newton approximant is written as
\begin{equation}\label{Newman_approx}
R^{\rm{New}}(x) = x\frac{p(x)-p(-x)}{p(x)+p(-x)} =  \frac{x\big{[} \prod_{k=1}^{n-1}(x+\alpha^k)-\prod_{k=1}^{n-1}(-x+\alpha^k)\big{]}}{\prod_{k=1}^{n-1}(x+\alpha^k)+\prod_{k=1}^{n-1}(-x+\alpha^k)}.
\end{equation}
The resulting error on the positive interval written  $\gamma(x) = x - R^{\rm{New}}(x)$ and satisfies the following inequality on $[\epsilon,1]$:
\begin{equation}
\gamma(x) = x - R^{\rm{New}}(x) = 2x \frac{\frac{p(-x)}{p(x)}}{1+\frac{p(-x)}{p(x)}} = 2x \frac{f(x)}{1+f(x)} \leqslant 2x \frac{f(x)}{1-f(x)} \leqslant  \frac{2f(x)}{1-f(x)},
\end{equation}
and using that $\vert f(x) \vert < e^{-\sqrt{n}}$ it follows that $\forall x \in [\epsilon,1]$, and $\ n \geqslant 5$
\begin{equation}
\gamma(x) < 3 e^{-\sqrt{n}}.
\end{equation}
Involving some additional reasoning, one can extend the above upper bound as shown in \cite{Ne64}, to the whole interval $[0,1]$, and because of symmetry reasons, to $[-1,1]$ as well.

\begin{remark}
	The Newman approximant $R^{\text{new}}(x)$ in (\ref{Newman_approx}) is indeed interpolatory. It can thus be obtained using the Loewner framework in \cite{ajm07}, with the following $n$ associated  interpolation pairs containing $2n-1$ interpolation points:
	\begin{equation}
	\{(0,0), \ \ (-\alpha,\alpha), \ \ (-\alpha^2,\alpha^2), \ \ \ldots, \ \ (-\alpha^{n-1},\alpha^{n-1})\}.
	\end{equation}
\end{remark}


\section{The main approximation methods}\label{sec:main}

\subsection{Linear dynamical systems}\label{sec:LTI}


A linear dynamical system $\Sigma$ is characterized by a set of differential equations provided below:
\begin{equation} \label{lin_sys} 
\Sigma:~~\bE \dot{\bx}(t)=\bA\bx(t)+\bB u(t),~~~~ 
y(t)=\bC\bx(t), 
\end{equation} 
where $\bE,\bA \in \IR^{n\times n}$, $\bB, \bC^T \in \IR^{n}$. The dimension of system $\Sigma$ is denoted with $n$ and it represents the length of the internal variable $\bx(t) \in \mathbb{R}^n$. The control input is denoted with $u(t) \in \mathbb{R}$, while the observed output is $y(t) \in \mathbb{R}$. The matrix $\bE$ is typically considered to be invertible and hence  assimilated into matrices $\bA$ and $\bB$ as follows: $\bar{\bA} = \bE^{-1} \bA, \bar{\bB} = \bE^{-1} \bB$. By applying the Laplace transform directly to the differential equation in (\ref{lin_sys}), one can write that
\begin{equation}
s \bE \bX(s) = \bA \bX(s) + \bB \bU(s) \Rightarrow \bX(s) = (s\bE-\bA)^{-1} \bB \bU(s) \Rightarrow \bY(s) =  \bC (s\bE-\bA)^{-1} \bB \bU(s)
\end{equation}
Hence, the transfer function is the ratio of the output $\bY(s)$ to the input $\bU(s)$:
\begin{equation}\label{TF_lin}
H(s) = \frac{\bY(s)}{\bU(s)} = \bC (s\bE-\bA)^{-1} \bB.
\end{equation}
\noindent
{\bf Model reduction} generically refers to a class of methodologies that aim at, in the case of linear systems, constructing reduced-order linear systems of the form:
\begin{equation}\label{lin_sys_red}
\hat\Sigma:~
\hat\bE \dot{\hat\bx} (t)= {\hat\bA}{\hat\bx}(t)+{\hat\bB}\bu(t),~~
\hat\by(t)  =  \hat\bC\hat\bx (t)+\hat\bD\bu(t),
\end{equation}
where~ $\hat\bE,\,\hat\bA\in\IR^{r\times r}$, ~$\hat\bB\in\IR^{r\times m}$,~$\hat\bC\in\IR^{p\times r}$, ~$\hat\bD\in\IR^{p\times m}$. Such systems in (\ref{lin_sys_red}) are used to replace those in (\ref{lin_sys}) for numerically expensive tasks such as real-time simulation, optimal control, etc. The dimension of the reduced-order model $\hat{\Sigma}$ is denoted with $r$ and it is typically much lower than that of system $\Sigma$. 

For linear dynamical systems such as the one in (\ref{lin_sys}), the general idea is to start from data sets that contain frequency response measurements. The goal, as it will be made clear in the next section, is to find rational functions that (approximately) interpolate these measurements. More precisely, through a system theoretical viewpoint, we seek reduced-order linear dynamical system that models these measurements.

For more details on state-of-the-art MOR methods and an extensive historical account of the existing literature, we refer the reader to \cite{ACA05}. Extensions to nonlinear systems can be found \cite{baur14}, while extensions to parametric systems are presented in \cite{benner15}. For a recent account of the developments of MOR (algorithms and theory), see \cite{MORbook17}. The case of interpolation-based MOR methods is extensively and exhaustively treated in \cite{abg_book}.

\subsection{The Loewner framework}\label{sec:Loew}

The Loewner framework (LF) is a data-driven model identification and reduction tool introduced in \cite{ajm07}. It provides a solution through rational approximation by means of interpolation. For more details, see the extensive tutorial paper \cite{ALI17}, the Ph.D. dissertation in \cite{Io13} that addresses connections of the LF to Lagrange rational interpolation, the Ph.D. dissertation in \cite{Go17} that addresses extensions of the LF to some special classes of nonlinear dynamical systems, and the recent book chapter \cite{kga20}.

The set-up of the LF is as follows: we are given $N>0$ sampling points $\cT = \{\tau_1, \tau_2, \ldots,\tau_N\}$ together with $N$ function evaluations at these points denoted with $\{f_1, f_2, \ldots,f_N\}$ are given. The original set of measurement pairs is denoted with $\mathfrak{D} = \{(\tau_\ell; \  f_\ell) \vert \ell =1,\cdots,N\}$.

 Assume that $N$ is an even positive integer and let $k = N/2$. The given data $\mathfrak{D}$ is first partitioned into the following two disjoint subsets
\begin{equation}\label{data_Loew}
\mathfrak{D} = \mathfrak{L} \cup \mathfrak{R}, \ \text{where} \ \begin{cases} {\rm right \ data}: \mathfrak{L} = \{(\lambda_i; \  w_i) \vert i=1,\cdots,k
\}, \\
 {\rm left \ data}: \ \ \mathfrak{R} =
\{(\mu_j; \ v_j) \vert j=1,\cdots,k\}. \end{cases}
\end{equation} 
The elements of the subsets $\mathfrak{L}$ and $\mathfrak{R}$ are as follows (for $1 \leq i,j \leq k$):
\begin{enumerate}
 \item  $\lambda_i$, $\mu_j$ $\in\IC$ are the right, and, respectively left sample points.
	\item  $w_i$, $v_j$ $\in\IC$ are the right, and, respectively left sample values.
\end{enumerate}
The problem can be formulated as follows: find a rational function $R(s)$ of order $r \ll k$, such that the following interpolation conditions are (approximately) fulfilled:
\begin{equation} \label{interp_cond}
R(\lambda_i)=w_i,~~~R(\mu_j)=v_j.
\end{equation}

\begin{remark}
Note that in a typical data-driven MOR setup, it is considered that the sample values  $w_i$, $v_j$ $\in\IC$ represent evaluations of the transfer function $H(s)$ (as defined in (\ref{TF_lin})) corresponding to an underlying linear dynamical model as in (\ref {lin_sys}). More precisely, one would assume that $H(\lambda_i)  =w_i$ and $H(\mu_j)  =v_j$ for all $1 \leq i,j \leq k$.
\end{remark}

The next step is to arrange the data (\ref{data_Loew}) into matrix format. The Loewner matrix $\IL \in\IC^{k\times k}$ and the shifted Loewner matrix $\sIL \in\IC^{k\times k}$ are defined as follows
\begin{equation} \label{Loew_mat}
\IL=\left[\begin{array}{ccc}
\frac{v_1-w_1}{\mu_1-\lambda_1} & \cdots &
\frac{v_1-w_k}{\mu_1-\lambda_k} \\
\vdots & \ddots & \vdots \\
\frac{v_k-w_1}{\mu_k-\lambda_1} & \cdots &
\frac{v_k-w_k}{\mu_k-\lambda_k} \\
\end{array}\right], \ \ 
\sIL=\left[\begin{array}{ccc}
\frac{\mu_1v_1-w_1\lambda_1}{\mu_1-\lambda_1} & \cdots &
\frac{\mu_1v_1-w_k \lambda_k}{\mu_1-\lambda_k} \\
\vdots & \ddots & \vdots \\
\frac{\mu_kv_k-w_1 \lambda_1}{\mu_k-\lambda_1} & \cdots &
\frac{\mu_kv_k-w_k\lambda_k}{\mu_k-\lambda_k} \\
\end{array}\right],
\end{equation}
while the data vectors $\bV, \bW^T \in \IR^k$ are introduced as
\begin{equation} \label{VW_vec}
\bV=\left[\!\begin{array}{cccc} v_1 \  \ v_2 \  \ \ldots \  \ v_k \end{array}\!\right]^T, \ \ \ \bW=[w_1~~w_2~~\cdots~~w_k].
\end{equation}
Additionally, introduce the diagonal matrices $\bLa, \bM$ containing the sample points as
\begin{align}\label{LAM_M_mat}
\bLa = \text{diag}(\left[ \begin{array}{ccc} \lambda_1 & \ldots & \lambda_k  \end{array} \right]) \in \mathbb{R}^{k \times k}, \ \ \ \bM = \text{diag}(\left[ \begin{array}{ccc} \mu_1 & \ldots & \mu_k  \end{array} \right]) \in \mathbb{R}^{k \times k},
\end{align}	
and the vectors of ones $\bR = \bL^T = [1~~1~~\cdots~~1] \in \mathbb{C}^{1 \times k}$.

The following relations between the Loewner and the shifted Loewner matrices were proven to hold (see for example in \cite{ajm07}): 
\begin{equation}\label{Loew_eq}
 \sIL =  \IL \bLa + \bV \bR, \ \ \ \ 
\sIL =  \bM \IL + \bL \bW.
\end{equation}
It was also shown in \cite{ajm07} that the Loewner matrices satisfy the following Sylvester equations:
\begin{equation}\label{sylv_eq}
\bM \IL - \IL \bLa = \bV \bR - \bL \bW, \ \ \
\bM \sIL - \sIL \bLa = \bM \bV \bR - \bL \bW \bLambda.
\end{equation}	
where $\bR = \bL^T = \left[ \begin{array}{ccc} 1 & \ldots & 1  \end{array} \right] \in \mathbb{R}^{1 \times k}$.

\begin{lemma}
	If the matrix pencil $(\sIL,\,\IL)$ is regular, then one can directly recover the system matrices and the interpolation function. The Loewner model is composed of
	\begin{equation}
	\bE=-\IL,~~ \bA=-\sIL,~~ \bB=\bV,~~ \bC=\bW, \quad \rightarrow \quad \bH(s)=\bW(\sIL-s\IL)^{-1}\bV.
	\end{equation}
\end{lemma}
Introduce $\cR$ as the generalized controllability matrix in terms of the right sampling points $\{\lambda_1,\ldots,\lambda_k\}$ and matrices $\bE,\bA,\bB$. Additionally, let $\cO$  be the generalized observability matrix, written in terms of the left sampling points $\{\mu_1,\ldots,\mu_k\}$ and matrices $\bE,\bA,\bC$.
\begin{equation} \label{defOR}
\cO = 	\left[\begin{array}{c}
\bC(\mu_1\bE-\bA)^{-1}\\
\vdots\\
\bC(\mu_k\bE-\bA)^{-1}
\end{array}
\right], \ \ \cC = \left[
\begin{array}{ccc}
(\lambda_1\bE-\bA)^{-1}\bB & \cdots & (\lambda_k\bE-\bA)^{-1}\bB
\end{array}\right].
\end{equation}
Then, it follows that the data matrices can be factorized as follows
\begin{equation}
\IL = -\cO \bE \cC, \ \ \sIL = - \cO \bA \cC, \ \bV=\bC \cC, \ \bW=\cO\bB. 
\end{equation}

In practical applications, the pencil $(\sIL,\,\IL)$ is often singular. In these cases, perform a rank revealing singular value decomposition (SVD) of the Loewner matrix $\IL$.
By setting ${\rm rank}\,\IL =k$, write ($\text{where} \ \ \bX_r = \bX(:,1:r), \ \bY_r = \bY(:,1:r)$, and $\bS_r = \bS(1:r,1:r)$)
\begin{equation}
~\IL =\bX \bS \bY^* \approx \bX_r\bS_r \bY_r^*, \ \ \ \text{with} \ \ \bX_r, \bY_r \in \IC^{k \times r}, \  \bS_r \in \IC^{r \times r}.
\end{equation}	
\begin{lemma}			
	The system matrices corresponding to a compressed /projected Loewner model for which the transfer function approximately matches the conditions in (\ref{interp_cond}), can be computed as
	\begin{equation} \label{Loew_red_mat}
	\hat{\bE}_{\rm Loew} = -\bX_r^*\IL \bY_r, \ \  \hat{\bA}_{\rm Loew} = -\bX_r^*\sIL \bY_r, \ \
	\hat{\bB}_{\rm Loew} = \bX_r^*\bV, \ \  \hat{\bC}_{\rm Loew} = \bW \bY_r.
	\end{equation}
	The rational function corresponding to the reduced-order Loewner model can be computed in the following way
	\begin{align}\label{loew_fct}
	\begin{split}
	R_{\rm Loew}(x) = \hat{\bC}_{\rm Loew} (x \hat{\bE}_{\rm Loew}-\hat{\bA}_{\rm Loew})^{-1} \hat{\bB}_{\rm Loew} =  \bW \bY_r \Big{(} \bX_r^* (\sIL- x \IL) \bY_r \Big{)}^{-1} \bX_r^* \bV.
		\end{split}
	\end{align}
\end{lemma}

\begin{remark}
Provided that the reduced-order matrix $\hat{\bE}$ in (\ref{Loew_red_mat}) is invertible, one can incorporate it into the other system matrices. In this way, we put together an equivalent model $(\bI_r, \tilde{\bA}_{\rm Loew},\tilde{\bB}_{\rm Loew},\tilde{\bC}_{\rm Loew}) \equiv (\hat{\bE}_{\rm Loew}, \hat{\bA}_{\rm Loew},\hat{\bB}_{\rm Loew},\hat{\bC}_{\rm Loew}) $ in "standard" representation:
	\begin{align*} \label{Loew_red_noE}
	\tilde{\bA}_{\rm Loew} &= \hat{\bE}_{\rm Loew}^{-1} \hat{\bA}_{\rm Loew} = (\bX_r^*\IL \bY_r)^{-1}  (\bX_r^*\sIL \bY_r), \\ \tilde{\bB}_{\rm Loew} &=- \hat{\bE}_{\rm Loew}^{-1} \hat{\bB}_{\rm Loew} = - (\bX_r^*\IL \bY_r)^{-1}  (\bX_r^* \bV), \  \tilde{\bC}_{\rm Loew} = \hat{\bC}_{\rm Loew} = \bW \bY_r.
	\end{align*}
\end{remark}

\subsubsection{Extension to odd number of measurements}

The Loewner framework does not necessarily impose a data set $\mathfrak{D}$ in (\ref{data_Loew}) with an even number of  measurements. Moreover, the Loewner framework was recently extended in \cite{aca16} to cope with transfer functions of rectangular systems.

By assuming that $N=2k+1$, i.e., we have one extra measurement pair denoted with $(\tau_{2k+1},f_{2k+1})$. Assume that this pair is to be included in the right data set in (\ref{data_Loew}), i.e. let $\lambda_{k+1} = \tau_{2k+1}$ and $w_{k+1} = f_{2k+1}$. The new updated Loewner matrices $\breve{\IL}$, and $\breve{\sIL}$ can be written in terms of the original Loewner matrices in (\ref{Loew_mat}) as follows
\begin{align}
\begin{split}
\breve{\IL} = \left[ \IL \ \  \bz_1 \right] \in \mathbb{C}^{k \times (k+1)}, \ \ \bz_1 \in \mathbb{C}^k, \ \ \bz_1(j) = \frac{v_j-w_{k+1}}{\mu_j-\lambda_{k+1}}, \\
\breve{\sIL} = \left[ \IL \ \  \bz_2 \right] \in \mathbb{C}^{k \times (k+1)}, \ \ \bz_2 \in \mathbb{C}^k, \ \ \bz_2(j) = \frac{\mu_j v_j-\lambda_{k+1} w_{k+1}}{\mu_j-\lambda_{k+1}}.
\end{split}
\end{align}
Additionally, the following vectors are given as:
\begin{equation}
\breve{\bV} = \bV, \ \ \ \ \breve{\bW} = \left[ \bW \ \ w_{k+1} \right] \in \mathbb{C}^{1 \times (k+1)}, \ \ \breve{\bR} = [\bR~~1] \in \mathbb{C}^{1 \times (k+1)}, \  \breve{\bLambda} = \left[ \begin{matrix}
\bLambda & \\ & \lambda_{k+1}
\end{matrix} \right], \ \breve{\bM} = \bM.
\end{equation}

\begin{remark}
	Note that the procedure discussed in this section will further be used in Section\;\ref{sec:numerics} which includes the numerical results. There, the origin point $0$ will be added to the set of $2k$ points chosen inside the interval $\mathcal{I} = [-b,-a] \cup [a,b]$, where $0 < a < b < 1$. The reason for adding $0$ is to decrease the overall approximation error.
\end{remark}
	
\subsubsection{Extension to identical left and right data sets}

The Loewner framework does not necessarily impose disjoint left and right data subsets in (\ref{data_Loew}). As explained in the original paper \cite{ajm07}, one can modify the classical Loewner framework by including information of the derivative corresponding to the underlying "to-be-approximated" function $H(s)$.

Here we show how to treat the case of identical subsets, i.e., $\mu_i = \lambda_i, \ \ \forall 1 \leq i \leq k$ and $v_i = w_i, \ \ \forall 1 \leq i \leq k$. We adapt the  definition of the Loewner matrices in (\ref{Loew_mat}) as was also proposed in \cite{bg12} (for the TF-IRKA iterative procedure)
\begin{align}
\underbar{$\IL$}(i,j) = \begin{cases}
\displaystyle \frac{w_i-w_j}{\lambda_i-\lambda_j}, \ \text{if} \ i \neq j \\[4mm]
\displaystyle \frac{dH(s)}{ds}\Big{\vert} (s = \lambda_i), \ \text{if} \ i = j.
\end{cases} \ \ \text{and} \ \ \underbar{$\sIL$}(i,j) = \begin{cases}
\displaystyle \frac{\lambda_i w_i-\lambda_j w_j}{\lambda_i-\lambda_j}, \ \text{if} \ i \neq j \\[4mm]
\displaystyle \frac{d[sH(s)]}{ds}\Big{\vert} (s = \lambda_i), \ \text{if} \ i = j.
\end{cases}
\end{align}
Additionally $\underbar{$\bW$}(i,j) = \underbar{$\bV$}^T(i,j) = [w_1~~w_2~~\cdots~~w_k] \in \mathbb{C}^{1 \times k}$.

\begin{remark}
	Note that the procedure discussed in this section will be used in Section\;\ref{sec:numerics} which includes the numerical results. There, the special case of data partitioning that will be referred to as  "Loewner same", will be based on this procedure.
\end{remark}

\subsection{The AAA algorithm}\label{sec:AAA}

The AAA (Adaptive Antoulas Anderson) algorithm proposed in \cite{nst18} represents an adaptive extension of the method originally introduced by A.C. Antoulas and B.D.O. Anderson in \cite{aa86}.

The AAA algorithm is a data-driven method that requires as an input only evaluations of a univariate function $F(s)$ on a set of points denoted with $\cT$. Extensions of the algorithm to multivariate functions have been recently proposed in \cite{rg20}.

Assume as given $N$ sampling points $\cT = \{\tau_1, \tau_2, \ldots,\tau_N\}$ together with $N$ function evaluations at these points denoted with $\{f_1, f_2, \ldots,f_N\}$, i.e. $F(\tau_i) = f_i$.

The AAA algorithm is based on an iteration. Assume we are at step $\ell \leq 1$. The order $(\ell,\ell)$ rational approximant $R_{\rm AAA}(x)$ that is constructed through the algorithm in \cite{nst18} is first written in the following barycentric format:
\begin{equation}\label{AAA_approx}
R_{\rm AAA}(x) = \frac{ \sum_{k=0}^\ell \frac{ \alpha_k w_k }{x - \nu_k }}{ \sum_{k=0}^\ell \frac{\alpha_k }{x- \nu_k}},
\end{equation}
where $\{\nu_0,\ldots,\nu_\ell\}$ are the interpolation points (selected from the set $\cT$), $\{w_1,\ldots,w_\ell\}$ are the function evaluations, and $\{\alpha_0,\ldots,\alpha_\ell\}$ are the weights.
The method enforces interpolation at the $\ell$ points $\{\lambda_0,\ldots,\lambda_\ell\}$, i.e. $R_{\rm AAA}(\nu_k) = w_k = F(\nu_k), \ \forall 0 \leqslant k \leqslant \ell$.

One needs to determine suitable weights $\alpha_1,\alpha_2,\ldots,\alpha_n$ by formulating the approximation problem as follows:
\begin{align}
\begin{split}
\min\limits_{\alpha_1,\ldots,\alpha_\ell} \sum_{i=1}^N (R_{\rm AAA}(\tau_i) - F(\tau_i) )^2& \Leftrightarrow \min\limits_{\alpha_1,\ldots,\alpha_\ell} \sum_{i=1}^{N-\ell} \Big{(} \frac{ \sum_{k=0}^\ell \frac{ \alpha_k w_k }{\tau_i- \nu_k }}{\sum_{k=0}^\ell \frac{\alpha_k }{\tau_i- \nu_k}} - f_i \Big{)}^2 \Leftrightarrow \
 \min\limits_{\alpha_1,\ldots,\alpha_\ell} \sum_{i=1}^{N-\ell}  \Big{(} \frac{\sum_{k=0}^\ell \frac{\alpha_k(w_k-f_i)}{\tau_i- \lambda_k}}{\sum_{k=0}^\ell \frac{\alpha_k }{\tau_i- \nu_k}} \Big{)}^2.
\end{split}
\end{align}
Instead of solving the more general problem above (which is nonlinear in variables $\alpha_1, \alpha_2, \ldots, \alpha_\ell$) one solves a relaxed problem by means of linearization. This procedure leads to a least squares minimization problem, i.e.,
\begin{equation}
\min\limits_{\alpha_1,\ldots,\alpha_\ell} \sum_{i=1}^{N-\ell} \Big{(} \sum_{k=0}^\ell \frac{\alpha_k(f_i-w_k)}{\tau_i- \nu_k} \Big{)}^2.
\end{equation} 
Introduce the Loewner matrix $\IL \in \mathbb{C}^{(N-\ell-1) \times (\ell+1)}$ as $\IL(i,k) = \frac{f_i-w_k}{\tau_i- \nu_k}, \forall 1 \leqslant i \leq N-\ell-1$ and $1 \leq k \leq \ell+1$. Additionally, let $\balpha \in \mathbb{C}^{\ell+1}$ with $\balpha(k) = \alpha_k$. The minimization problem is written as
\begin{equation}
\min\limits_{\balpha} \Vert \IL \balpha \Vert_2^2.
\end{equation}
The solution is given by the $\ell+1$th right singular vector of the matrix $\IL$. The next interpolation point $\nu_{l+1}$ is selected by means of a Greedy approach, i.e., the location in $\Omega  \setminus \{\nu_0,\cdots,\nu_{\ell}\}$ where the error $\epsilon(s) = \vert F(s) - R(s) \vert, \ \ s \in \Omega$ is maximal.

Assume that at the end of the iteration, the rational approximant is of order $(r,r)$. Hence, one can also explicitly write this function $R_{\rm AAA}(x)$ as was done in the Loewner framework for the function  $R_{\rm Loew}(x)$ in (\ref{loew_fct}), i.e.,
\begin{equation}\label{AAA_approx2}
R_{\rm AAA}(x) = \hat{\bC}_{\rm AAA} (x \hat{\bE}_{\rm AAA}-\hat{\bA}_{\rm AAA})^{-1} \hat{\bB}_{\rm AAA},
\end{equation}
where the matrices are of dimension $r+2$ and are given explicitly in terms of the selected support points  $\{\lambda_0,\ldots,\lambda_r\}$, the function evaluations $\{w_1,\ldots,w_r\}$, and the weights $\{\alpha_0,\ldots,\alpha_r\}$. 
\begin{align}\label{AAA_matrices}
\begin{split}
 \hat{\bA}_{\rm AAA} &= \text{diag}(\left[\nu_0,\nu_1,\ldots,\nu_r,1\right]) -  \hat{\bB}_{\rm AAA}  \bfe_{r+2}^T, \ \ \hat{\bB}_{\rm AAA} = \left[ \begin{matrix}
\alpha_0 & \cdots & \alpha_r & 1 
\end{matrix}  \right]^T, \\ \hat{\bE}_{\rm AAA} &= \text{diag}(\left[1,1,\ldots,1,0\right]), \ \ \ \hat{\bC}_{\rm AAA} = \left[ \begin{matrix}
w_0 & \cdots & w_r & 0 
\end{matrix}  \right], \ \ \ \bfe_{r+2} = \left[ \begin{matrix}
1 & 1 & \cdots & 1
\end{matrix}\right]^T.
\end{split}
\end{align}
Note that for the numerical examples presented in Section\;\ref{sec:numerics}, we will be using the AAA implementation available in the Chebfun toolbox \cite{chebfun}.

\subsection{The minimax algorithm}\label{sec:minimax}

Computing the best polynomial approximation is a known problem that was solved many decades ago by Evgeny Yakovlevich Remez  who proposed an iterative algorithm (see \cite{Re34}). A robust implementation of the linear version of the Remez algorithm was recently proposed in \cite{pt09} (the code is available in the Chebfun toolbox \cite{chebfun}). Now switching to rational approximation, it is to be noted that there exist algorithms which address best rational approximation: the rational Remez algorithm and differential correction (DC) algorithm. In \cite{fntb18}, a new robust/optimized algorithm is proposed (based on adaptive construction of rational approximants in barycentric representations).

The problem that is addressed in the contribution \cite{fntb18} is to approximate functions $f \in \cC([a,b])$ using type $(m,n)$ rational approximations with real coefficients from the set $\cR_{m,n}$.

Given the function $f$ and integers $m,n$, the problem is formulated as follows: find $r \in \cR_{m,n}$ so that the infinity norm $\Vert f - r \Vert_\infty = \max_{x \in [a,b]} \vert f(x) - r(x) \vert$ is minimized.
It is known that such a minimizer exists and is unique (Chapter 24 in  \cite{Tr13}).

Denote with $\overline{r} \in \cR_{m,n}$ the best (minimax) approximation of function $f$ (the solution to the above problem). It follows that there exists a sequence of ordered nodes where the function $f-r^*$ takes its global extreme value over the interval $[a,b]$ (with alternating signs). This sequence is denoted with $(x_0,x_1,\ldots,x_{m+n+1})$ with  $a \leqslant x_0 < x_1 < \cdots < x_{m+n+1} \leqslant b$, hence
\begin{equation}
f(x_k)-\overline{r}(x_k) = (-1)^{k+1} \Vert f - \overline{r} \Vert_{\infty}, \ \ \ k \in \{0,1,\ldots,m+n+1\}.
\end{equation}

Computing rational minimax approximations is a challenging task especially when there are singularities on or near the interval of approximation. In the recent contribution \cite{fntb18}, the authors propose different combinations of robust algorithms that are based on rational barycentric representations. One key feature is that the support points are chosen in an adaptive fashion as the approximant is computed (as for the AAA algorithm in \cite{nst18}).

In the numerical examples presented in this contribution we will be using the Chebfun minimax implementation available in the toolbox \cite{chebfun}. As discussed in \cite{fntb18} there, the method is based on using an AAA-Lawson method of iteratively re-weighted least-squares followed by a classical Remez algorithm to enforce generically quadratic convergence.

\begin{figure}[ht]
	\begin{center}
		\includegraphics[scale=1.2]{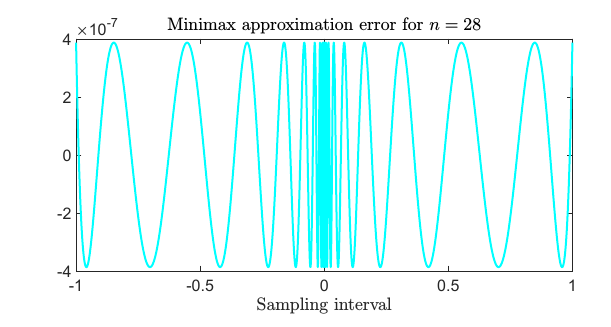}
	\end{center}
	\vspace{-3mm}
	\caption{Approximation curve produced by minimax on $[-1,1]$ - order $(28,28)$ approximant.}
	\label{fig01}
\end{figure}

In Fig.\;\ref{fig01}, we depict the magnitude of the approximation error for the best order (28,28) approximant (computed by means of the minimax algorithm in \cite{fntb18}). 

\begin{figure}[ht]
	\begin{center}
	\includegraphics[scale=0.4]{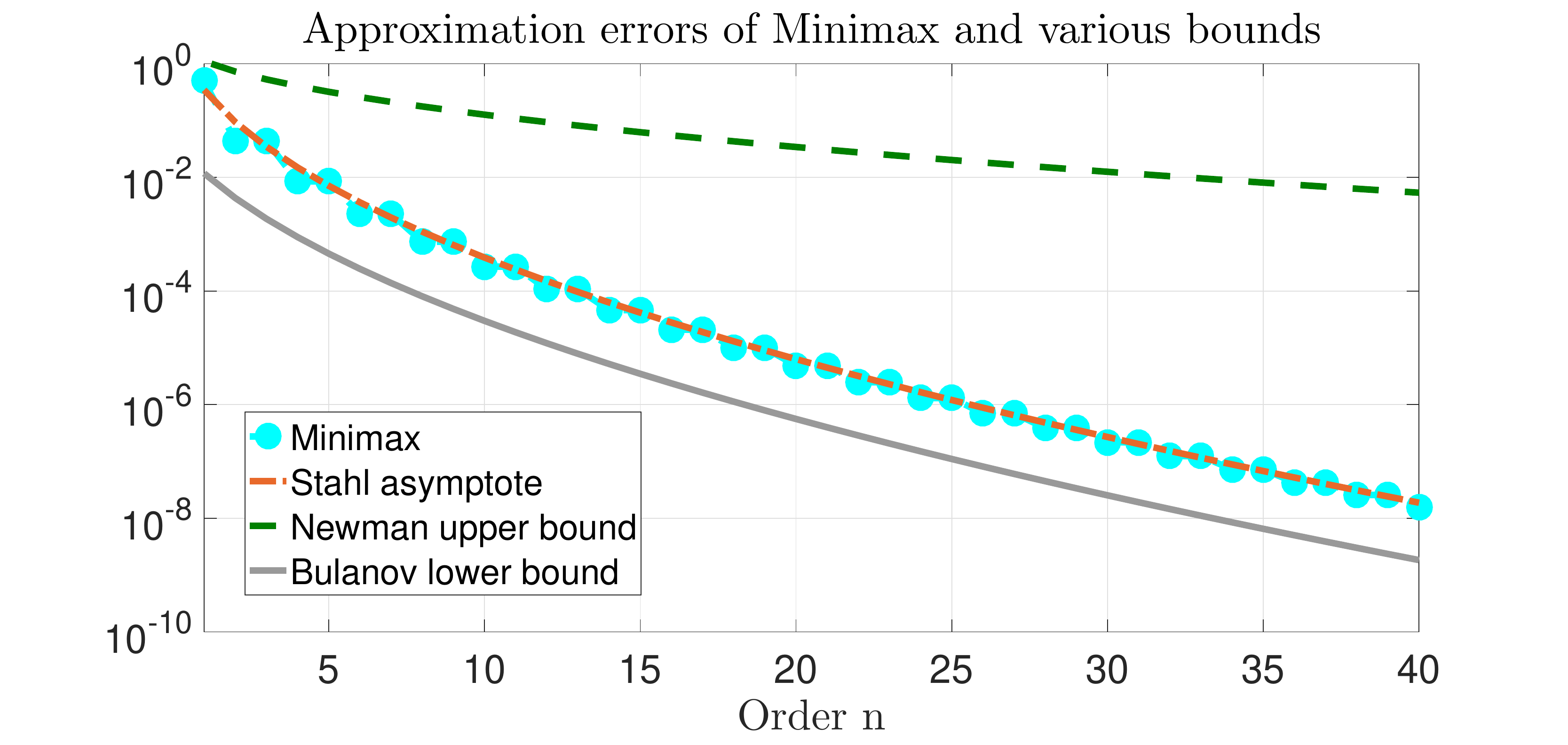}
	\end{center}
	\vspace{-3mm}
	\caption{Approximation errors produced by minimax on $[-1,1]$ for an order $(n,n)$ approximant with $1\leq n \leq 40$.}
	\label{fig1}
\end{figure}

Next, in Fig.\;\ref{fig1}, we depict the maximum approximation error in  for the best order (n,n) approximant computed by means of the minimax algorithm in \cite{fntb18}, where $1\leqslant n \leqslant 40$. Additionall, the bounds proposed by Newman, Bulanov and Stahl are shown on the same figure.

\subsection{Computation of the maximum error}\label{sec:error}

In this section we provide a brief description of the methodology used to compute the maximum error of a rational function $R(x)$ when approximating $\vert x \vert$ on $[-1,1]$. Such function $R(x)$ could be either computed with the Loewner framework introduced in Section\;\ref{sec:Loew}, with the AAA algorithm in Section\;\ref{sec:AAA}, or with the minimax algorithm in Section\;\ref{sec:minimax}. For the later, the maximum error is computed to be the absolute value of the $(r,r)$ rational minimax function evaluated at $0$.

Now, for the rational functions $R(x)$ computed with the Loewner and AAA approaches, we consider $(\hat{\bE},\hat{\bA},\hat{\bB}, \hat{\bC})$ to be a realization of the function $R(x)$. For the Loewner  framework, these matrices can be computed as in (\ref{Loew_red_mat}), while for the AAA algorithm, the matrices in (\ref{AAA_matrices}) are used instead.

Given a rational function $R(x) = \hat{\bC} (x \hat{\bE}-\hat{\bA})^{-1} \hat{\bB}$, we compute the maximum approximation error $\epsilon_R = \max_{x \in [-1,1]} \vert R(x) - |x|\vert$. Introduce also $\epsilon_R^+ = \max_{x \in (0,1)} \vert R(x) - |x|\vert$ and $\epsilon_R^- = \max_{x \in (-1,0)} \vert R(x) - |x|\vert$. Finally, for $\beta \in \{-1,0,1\}$, let $\epsilon_R^{\beta} =  \vert R(\beta) - |\beta|\vert$. Consequently, it follows that
\begin{equation}
\epsilon_R = \max(\epsilon_R^-,\epsilon_R^+,\epsilon_R^{-1},\epsilon_R^{1},\epsilon_R^0)
\end{equation}
The challenging part is computing $\epsilon_R^-$ and $\epsilon_R^+$, since the other 3 quantities are simple evaluations of $\vert R(x) - |x|\vert$ for  $x \in \{-1,0,1\}$. The main idea for computing $\epsilon_R^{+}$ and $\epsilon_R^{-}$ is as follows: determine all  extrema of the  function $R(x) - |x|$ by solving a generalized eigenvalue problem and select only the ones that are on the real axis in between $-1$ and $1$. 

We now go through the procedure for finding $\epsilon_R^+$ (finding $\epsilon_R^-$ can be derived i a similar manner). Finding all extrema of the  function $R(x) - x$ on $(0,1)$ is equivalent to finding the zeros of the rational function $R'(x) - 1$ on $(0,1)$. This can be formulated as finding the eigenvalues of the pair $(\overline{\bX},\overline{\bE})$, where the corresponding matrices can be found as follows
\begin{equation}
\overline{\bE} = \left[ \begin{matrix}
\hat{\bE} &  & \\  & \hat{\bE} & \\ & & 0
\end{matrix}  \right], \ \ \overline{\bA} = \left[ \begin{matrix}
\hat{\bA} & \hat{\bE} \\ \bfz & \hat{\bA}
\end{matrix}  \right], \ \  \overline{\bB} = \left[ \begin{matrix} \bfz \\ \hat{\bB} \end{matrix}  \right], \ \  \overline{\bC} = \left[ \begin{matrix} \hat{\bC} & \bfz \end{matrix}  \right],  \ \overline{\bD} = 1,\   \overline{\bX} = \left[ \begin{matrix}
\overline{\bA} & \overline{\bB} \\ \overline{\bC} & \overline{\bD}
\end{matrix}  \right].
\end{equation}
Finally, from the computed eigenvalues, select the one in $(0,1)$ at which the evaluation of $\vert R(x) - x \vert$ is maximal. It means that this is precisely the value we were looking for, i.e.,  $\epsilon_R^+$. Now, in order to compute  $\epsilon_R^-$, we repeat the above procedure by choosing instead $\overline{\bD} = -1$. 

\newpage

\section{Numerical experiments}\label{sec:numerics}

In this section we report on a variety of numerical experiments that were performed for approximating de absolute value function on the interval $[-1,1]$. Both the Loewner framework and the AAA algorithm are data-driven methods and hence they require measurements. We propose a couple of different types of sampling points and also different data partitioning techniques for the Loewner framework. Additionally, we provide comparisons to the classical rational approximation bounds derived in the previous century, which were also mentioned in Section\;\ref{sec:approx}.

\noindent
The sampling points are selected to follow five different distributions:
\begin{align}\label{sample_pts_types}
\begin{split}
& \bullet \text{\textbf{Linearly-spaced} \ (linspace) \ points;} \hspace{8mm} \bullet \text{\textbf{Logarithmically-spaced} \ (logspace) \ points;}  \\
& \bullet \text{\textbf{Chebyshev} \ points;} \hspace{25mm} \bullet \text{\textbf{Zolotarev}  \ points;}  \hspace{25mm} \bullet \text{\textbf{Newman}  \ points;}
\end{split}
\end{align}

\begin{figure}[ht]		\hspace*{-3mm} \includegraphics[scale=0.28]{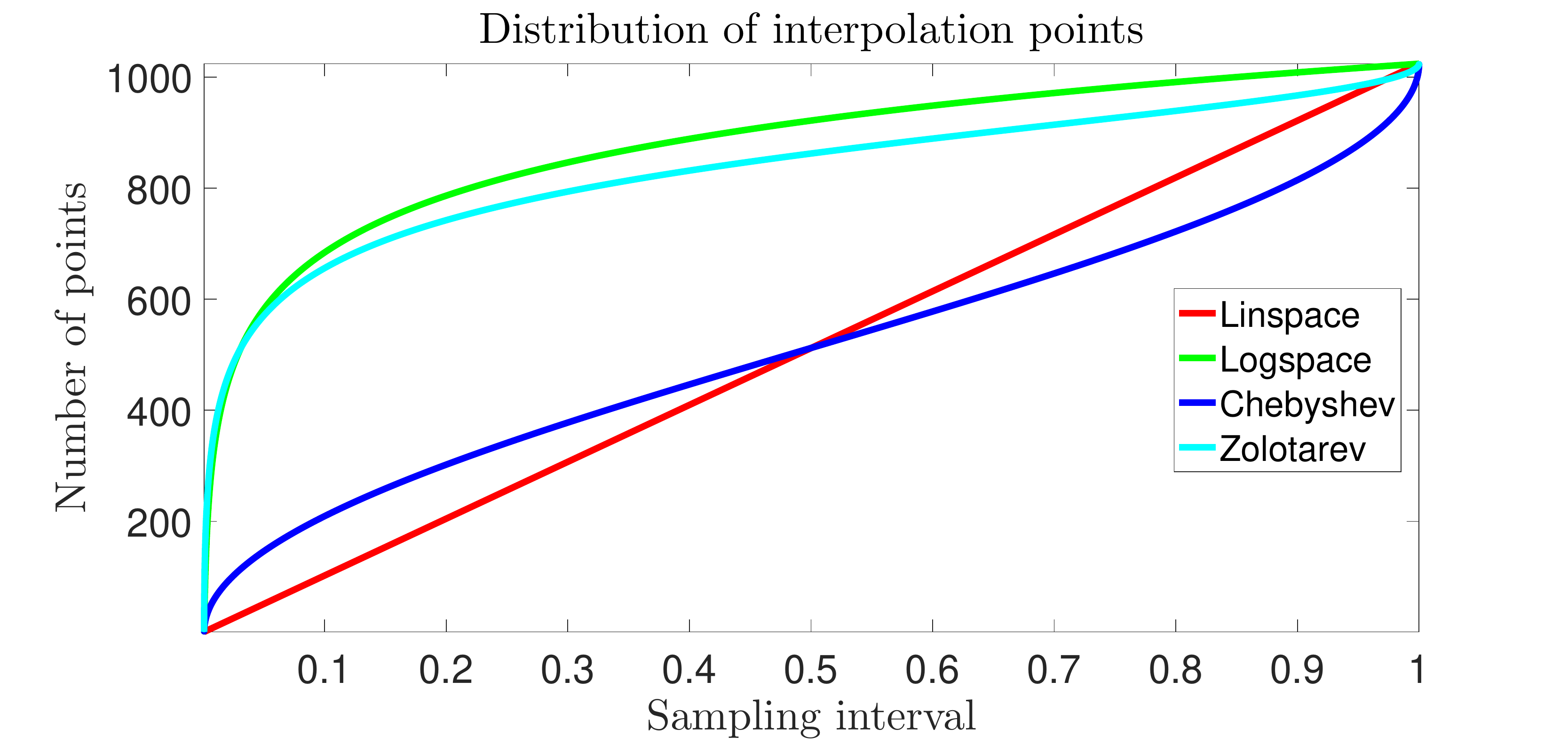}\hspace*{-9mm}	\includegraphics[scale=0.28]{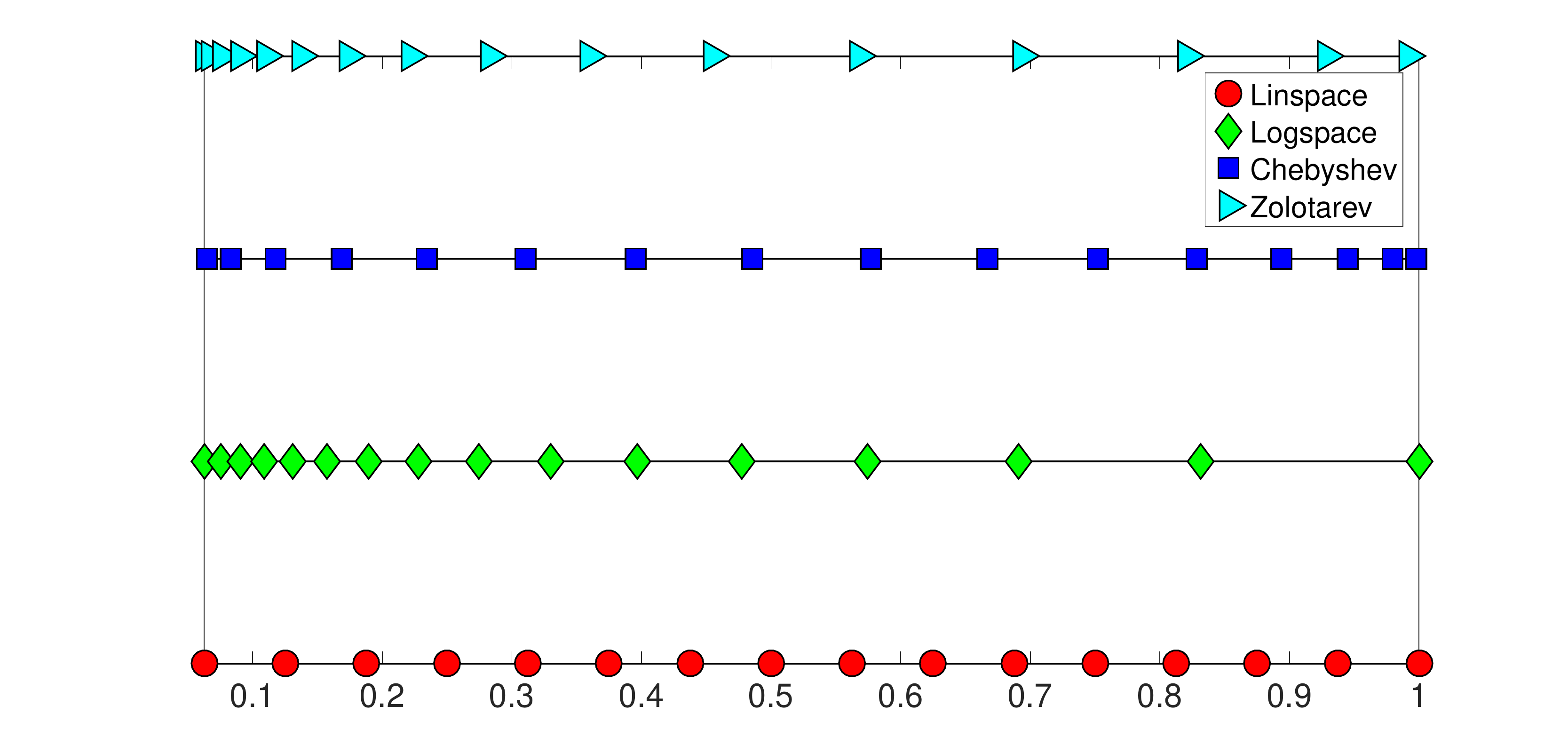}
	\vspace{-5mm}
	\caption{Four categories of sampling (interpolation points). Left: 1024 points (sampling interval vs number of points). Right: 16 points (on 4 parallel real axes).}
	\label{fig32}
	\vspace{-4mm}
\end{figure}
In Fig.\;\ref{fig32}, we depict four categories of sampling points with 2 examples (linspace, logspace, Chebyshev and Zolotarev). In the first one (left figure), 4 curves are plotted. The x axis represents the sampling interval $[0,1]$, while the y axis shows the number of points (in total 1024 points were chosen). In the second one (right figure), we depict 16 points for each of the four category with colorful symbols (the color of the curves coincide with those of the symbols and will remain consistent throughout the remaining of the manuscript).

Note that the partitioning of data into left and right data sets (see (\ref{data_Loew})) is the first step in the Loewner framework and a very crucial one indeed. In what follows, we use three data partition techniques for the Loewner framework:
\begin{align}\label{partition_types}
\begin{split}
&\bullet \text{\textbf{split} \ (disjoint \ intervals, \ e.g., \ negative \ and \ positive)}; \\
&\bullet \text{\textbf{alternating} \ (every \ second \ left \ and \ every \ second \ right)}; \\
&\bullet \text{\textbf{same} \ (the \ left \ and \ right \ points \ are \ the \ same)}.
\end{split}
\end{align}

\begin{figure}[ht]
	\begin{center}
		\includegraphics[scale=0.36]{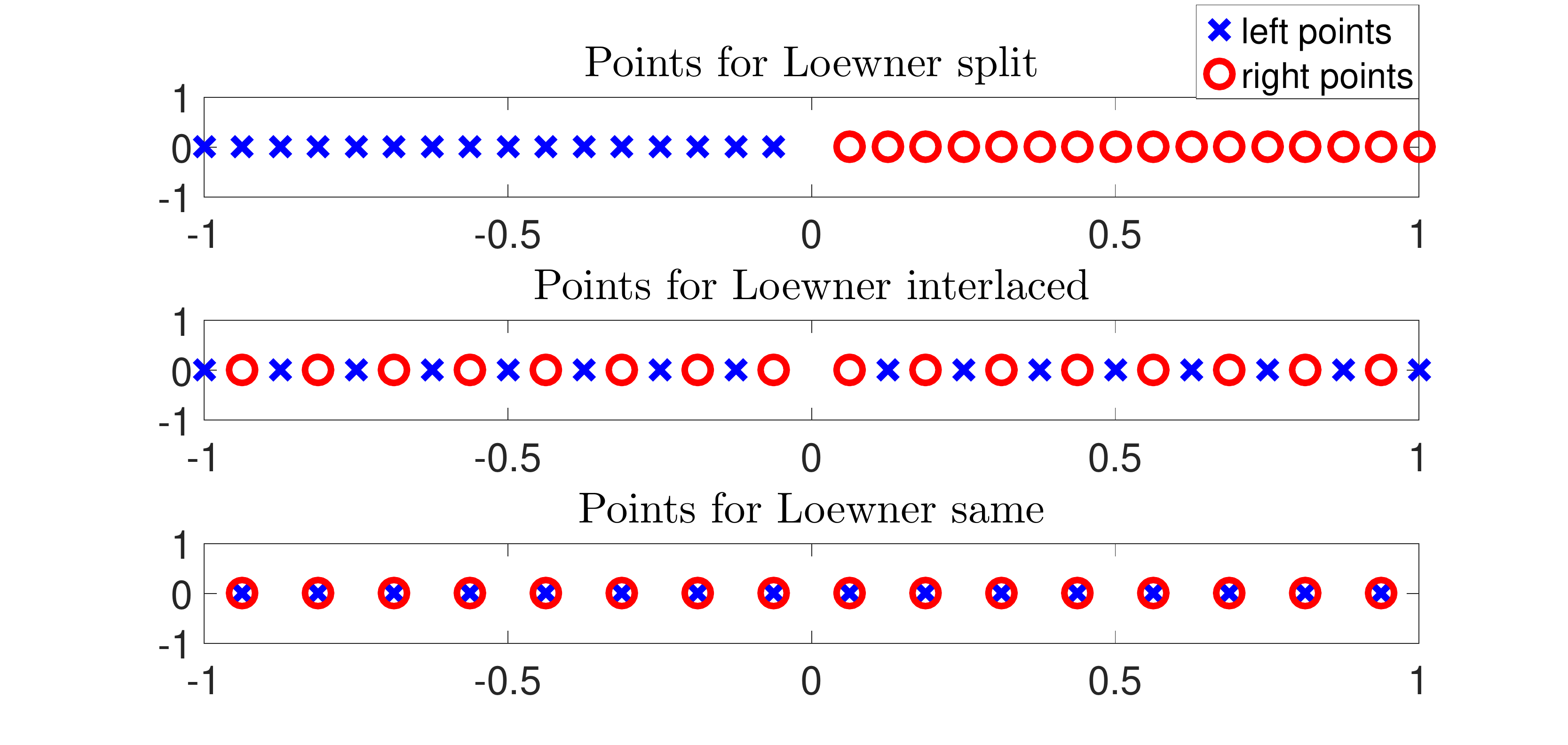}
	\end{center}
	\vspace{-3mm}
	\caption{The three partitioning options depicted for 16 linearly-spaced interpolation points.}
	\label{fig33}
\end{figure}

In Fig.\;\ref{fig33}, we depict the three partitioning options with an example of 16 points linearly spaced points in the interval $(a,b)$.
Let $a=2^{-10}$ and $b=1$ and choose $N=2048$ interpolation points inside $[-b,-a] \cup [a,b]$. More precisely, $k=1024$ are chosen inside the interval $[a,b]$ and the other half inside $[-b,-a]$.

In what follows, we apply the Loewner framework for the five types of sample points mentioned in (\ref{sample_pts_types}). For each of the different points selection, we apply three different partitioning schemes proposed in (\ref{partition_types}).

Note that the dimension of the Loewner matrices is $1024 \times 1024$ and the Loewner rational approximants are of order $(r-1,r)$. The AAA approximant is of order $(r,r)$. We will choose the reduction order to be $r=28$ for the following experiments.

\subsection{Linearly-spaced points}

For the first case, we consider that the interpolation points are linearly spaced within the interval $[a,b]$ for $a=2^{-10}$ and $b=1$. Let $n >0$ an integer and define positive "linspace" points $p^{\rm lin}_1 < p^{\rm lin}_2 < \cdots p^{\rm lin}_n$ as
\begin{equation}
p_1^{\rm lin} = a, \ \ p_2^{\rm lin} = a+\frac{b-a}{n-1}, \ \ \ldots \ \ ,  p_n^{\rm lin} = a+\frac{(n-1)(b-a)}{n-1} = b .
\end{equation}
Note that these points span a subinterval of the interval $(0,1)$. Additionally, for symmetry reasons, the negative points are also considered, i.e. $p_{-k}^{\rm lin} = - p_k^{\rm lin}$.

\begin{figure}[ht]		\includegraphics[scale=0.28]{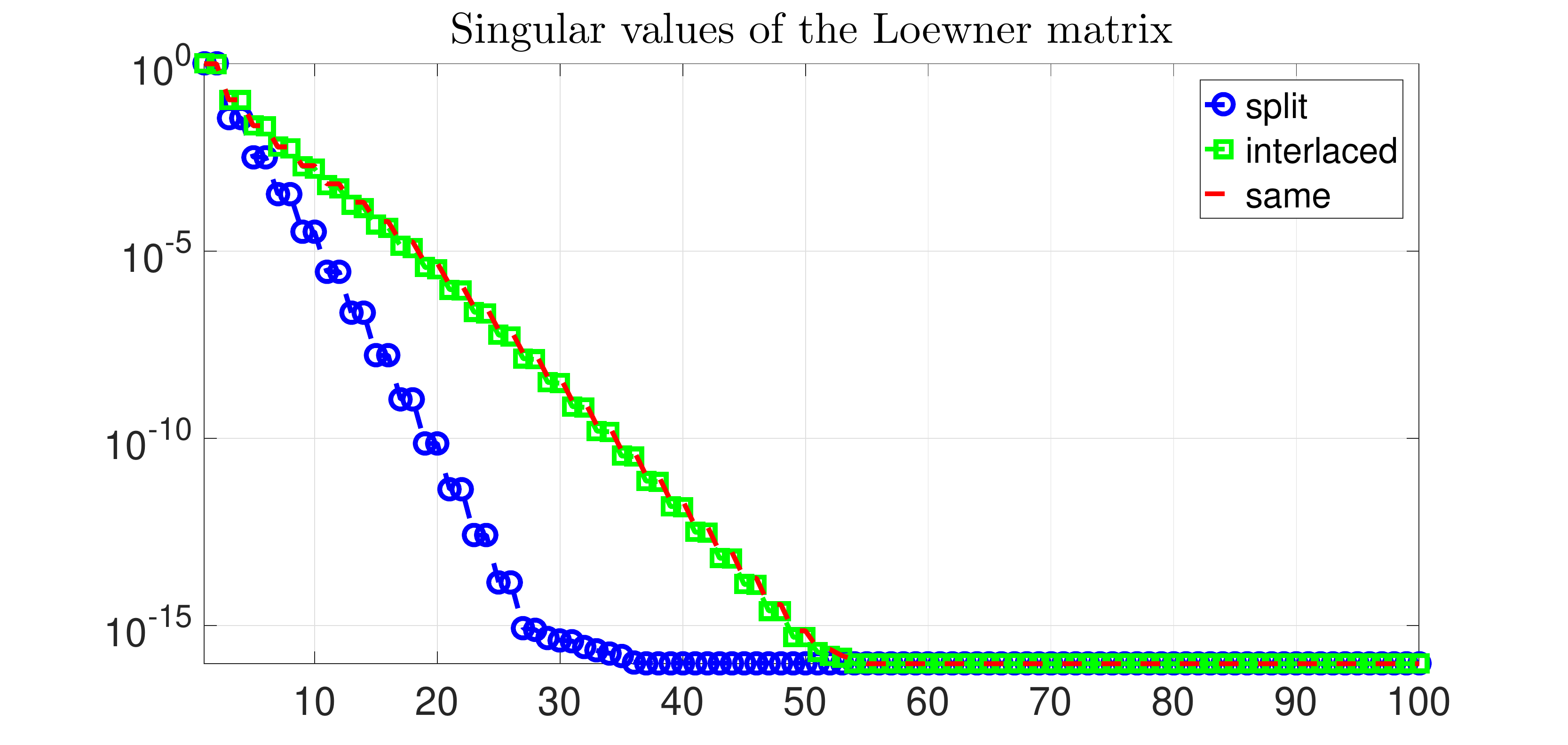}\hspace*{-8mm}	\includegraphics[scale=0.27]{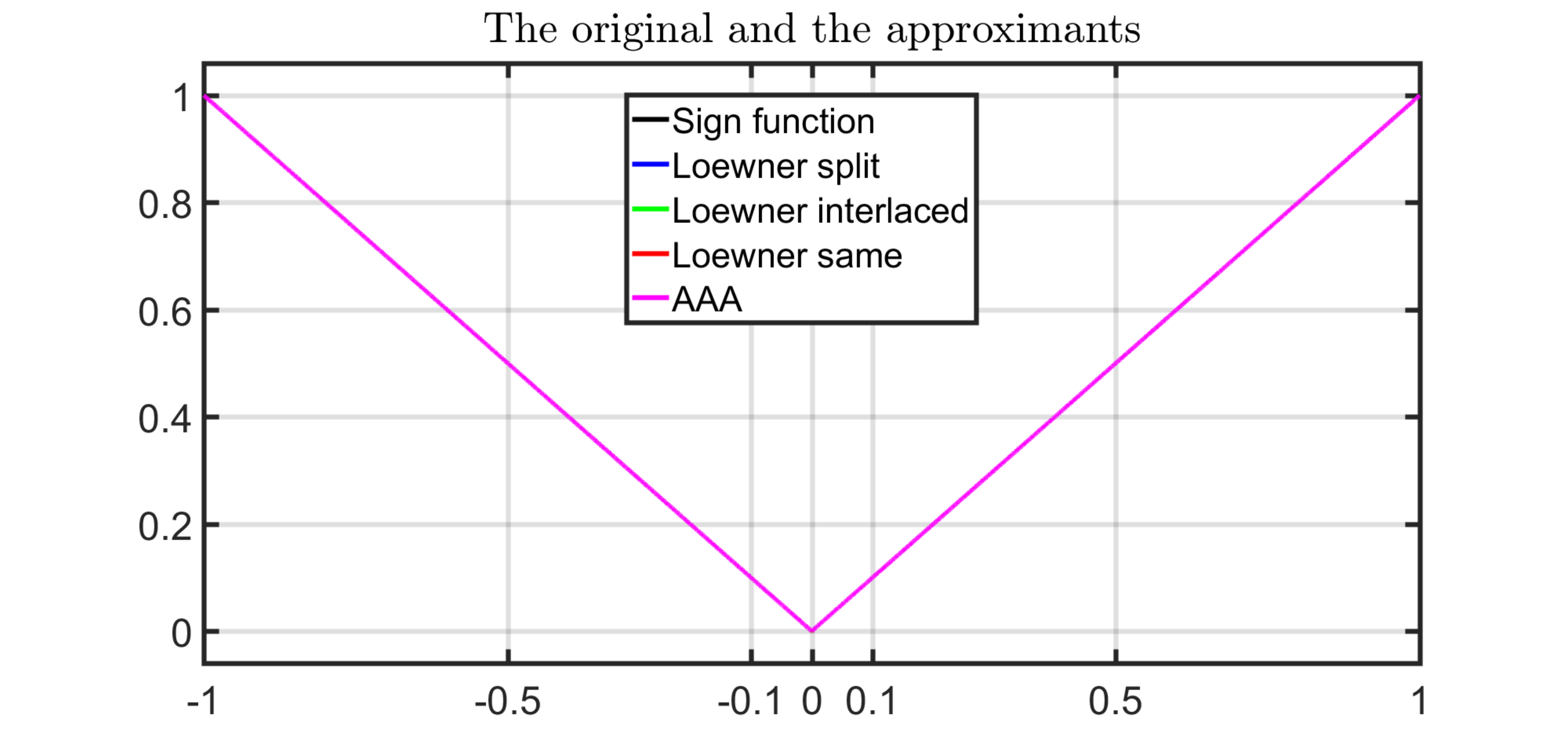}
	\vspace{-5mm}
	\caption{The first case (linspace points). Left: first 100 singular values of the 3 Loewner matrices. Right: the original function $\vert x \vert$ and the the 4 rational approximants.}
	\label{fig4}
	\vspace{-4mm}
\end{figure}

In Fig.\;\ref{fig4}, the singular value decay of the Loewner matrices computed by means of the three proposed partition techniques in (\ref{partition_types}) (left sub-figure). Moreover, the original absolute value function $\vert x \vert$ together with the 3 rational approximants computed with the Loewner framework (of order (27,28)) as well as the one of order (28,28) computed  by means of AAA, are depicted in the right sub-figure in Fig.\;\ref{fig4}.

\begin{figure}[ht]		
	\begin{center}
	\includegraphics[scale=0.94]{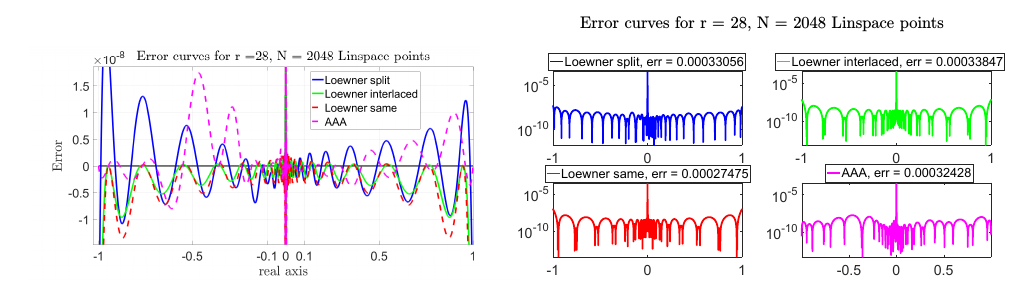}
	\end{center}
	\vspace{-5mm}
	\caption{The first case (linspace points). Left: the approximation errors for the 4 methods. Right: the magnitude of the errors for the 4 approximants.}
	\label{fig5}
	\vspace{-2mm}
\end{figure}

In Fig.\;\ref{fig5}, we depict the approximation error produced by the four rational functions that are used to approximate $\vert x \vert$. In the right sub-figure, the magnitude of the errors is shown, while in the left sub-figure all 4 error curves are depicted on the same plot. We observe that, in all four cases, there is a visible peak at 0.

\begin{figure}[ht]		
	\begin{center}
		\includegraphics[scale=1.3]{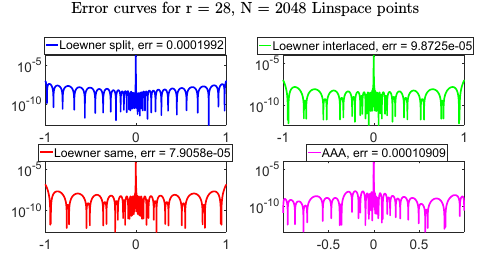}
	\end{center}
	\vspace{-5mm}
	\caption{The first case (linspace points). The magnitude of errors for the 4 approximants when 0 was added to the set of points.}
	\label{fig6}
	\vspace{-3mm}
\end{figure}

For the next experiment, we include $(0,0)$ in the set of $N = 2048$ pairs of measurements denoted with $\mathfrak{D}$. We apply again the Loewner framework for the three different partitioning options and compute the corresponding rational approximants (together with that computed using AAA). We depict the magnitude for all four error curves in Fig.\;\ref{fig6}.

\begin{figure}[ht]		
	\begin{center}
		\includegraphics[scale=0.36]{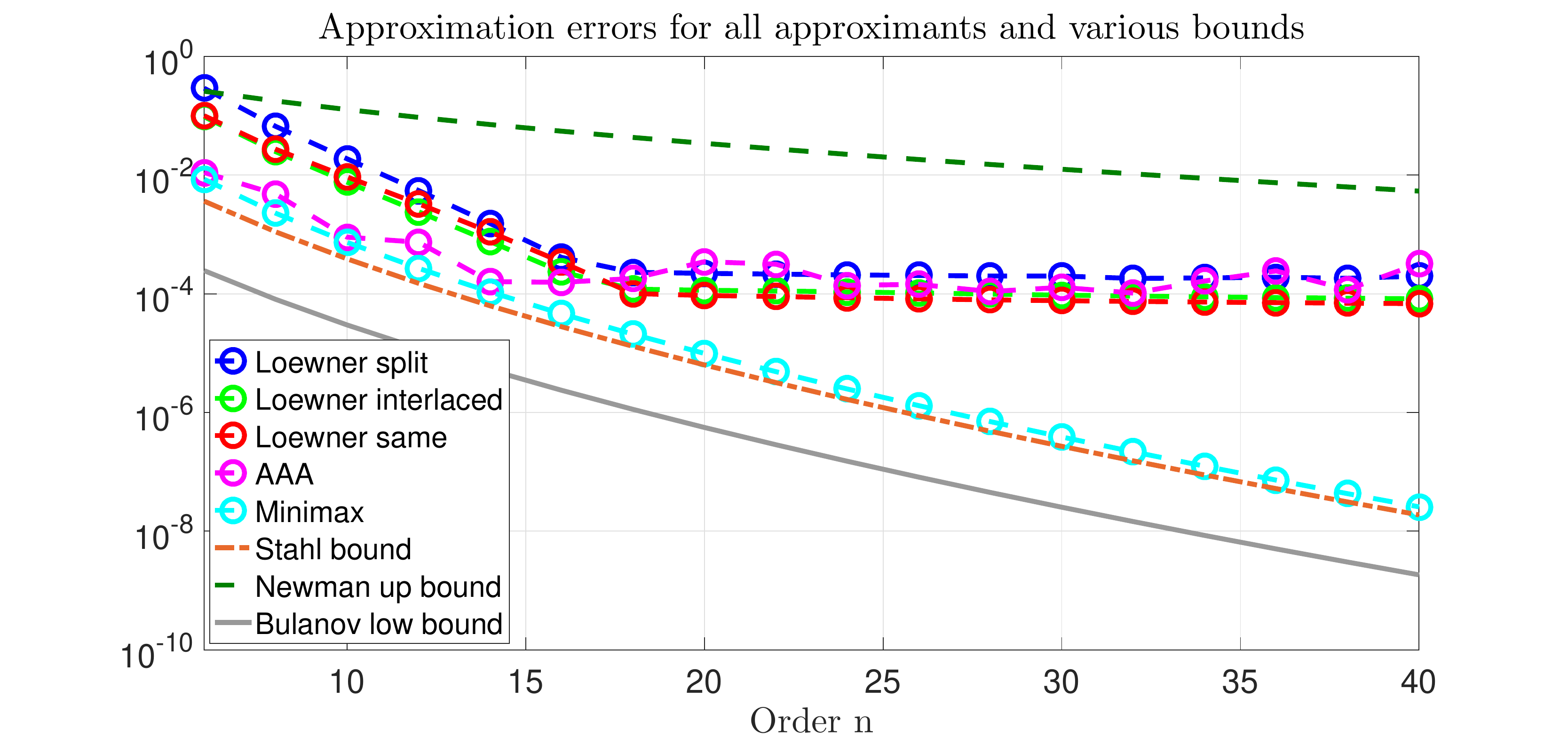}
	\end{center}
	\vspace{-5mm}
	\caption{The first case (linspace points). The maximum approximation errors for all methods and for different values of order $n$ (0 was added to the set of points)  + various bounds.}
	\label{fig7}
	\vspace{-3mm}
\end{figure}

Finally, by varying the reduction order from 6 to 40, we compute approximants by means of the four methods (note that 0 is again added here). We collect the maximum approximation errors for each method computed by means of the procedure presented in Section\;\ref{sec:iterate}. Additionally, we compute the Newman, Bulanov and Stahl bounds presented in Section\;\ref{sec:approx} (for each value). The results are depicted in Fig.\;\ref{fig7}.

Clearly, by adding the pair $(0,0)$ to the original set of measurement pairs $\mathfrak{D}$ of length $N=2048$ will improve the accuracy of the approximation. For example, for the Loewner same method, the approximation error decreased from $2.7575e-04$ to $7.9058e-05$.

Throughout the rest of the numerical experiments presented in this section, it will be assumed that 0 is included in the set of sample points.

\subsection{Logspace points}

In this section we consider that the interpolation points are logarithmically spaced within the interval $[-b,-a] \cup [a,b]$ for $a=2^{-10}$ and $b=1$. Let $\lambda_a = \log_{10}(a)$ and  $\lambda_b = \log_{10}(b)$. Let also $n >0$ an integer and define positive "logspace" points $p^{\rm log}_1 < p^{\rm log}_2 < \cdots p^{\rm log}_n$ as
\begin{equation}
p_1^{\rm log} = 10^{\lambda_a} = a, \ \ p_2^{\rm log} = 10^{\lambda_a +\frac{\lambda_b-\lambda_a}{n-1}}, \ \ \ldots \ \ ,  p_n^{\rm log} = 10^{\lambda_a +\frac{(n-1)(\lambda_b-\lambda_a)}{n-1}} = 10^{\lambda_b} =  b .
\end{equation}
Note that these points span a subinterval of the interval $(0,1)$. Additionally, for symmetry reasons, the negative points are also considered, i.e. $p_{-k}^{\rm log} = - p_k^{\rm log}$.

\begin{figure}[ht]	
		\includegraphics[scale=0.27]{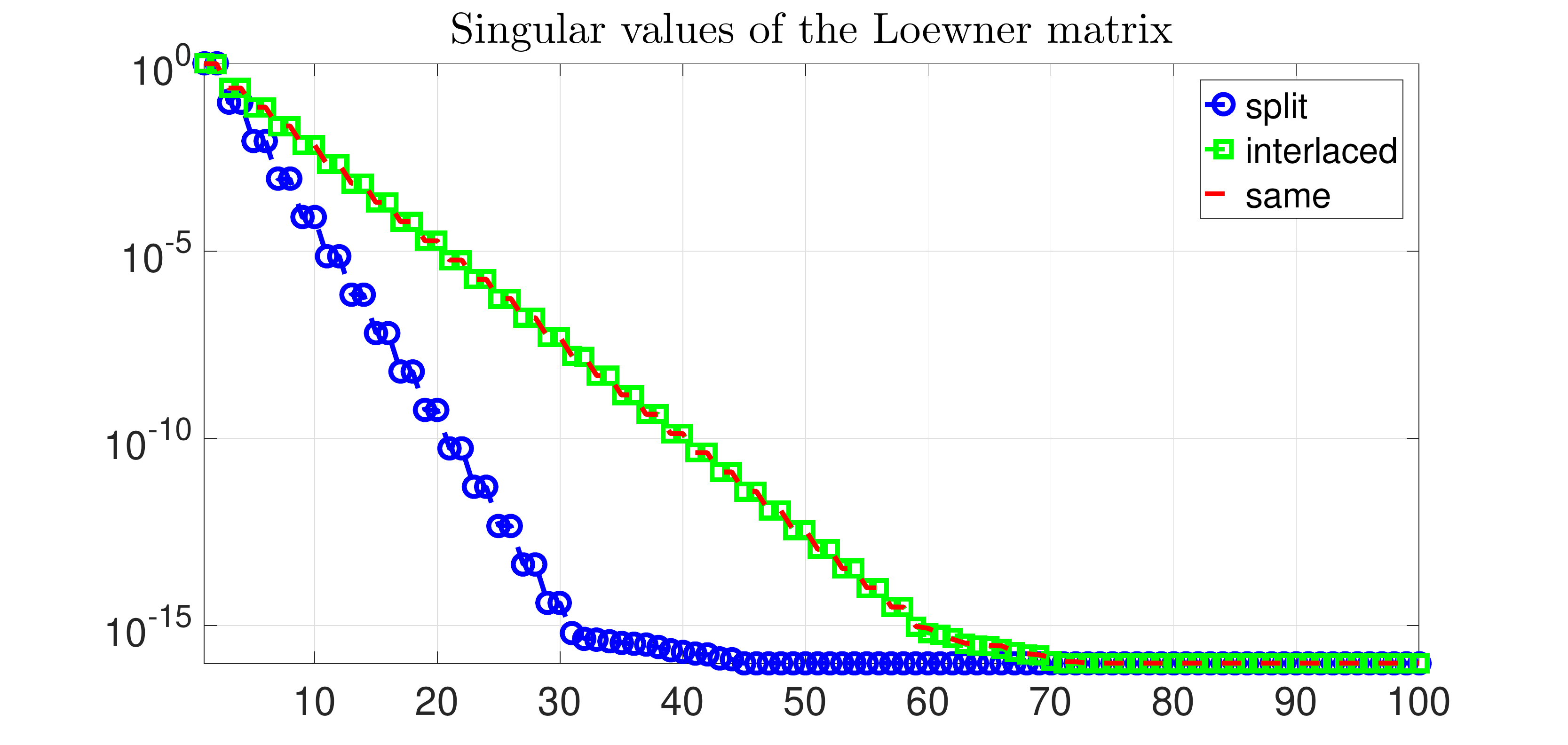}\hspace*{-8mm}	\includegraphics[scale=0.28]{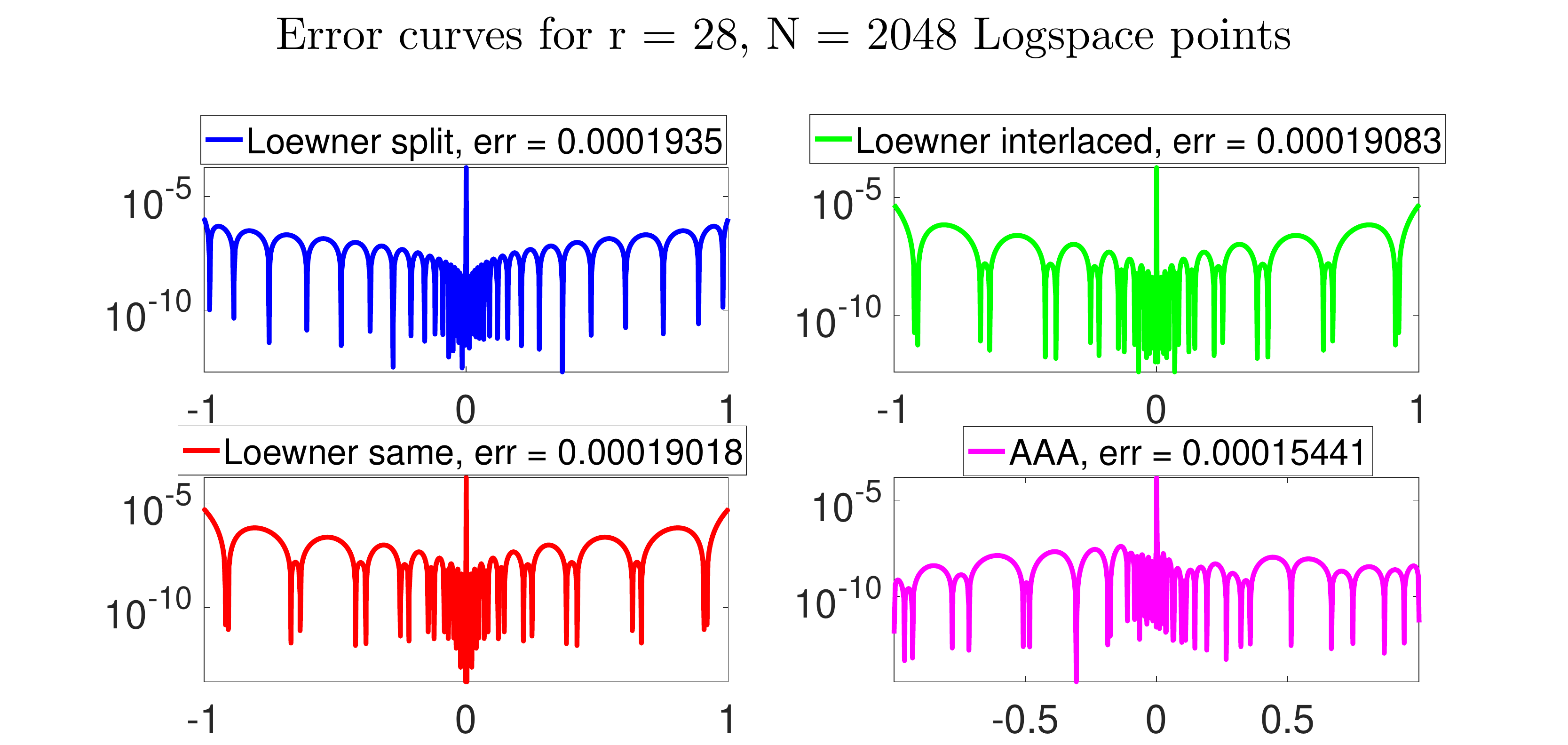}
	\vspace{-3mm}
	\caption{The second case (logspace points). Left: singular value decay of the Loewner matrices. Right: the magnitude of the error curves produced by order $r=28$ rational approximants.}
	\label{fig10}
	\vspace{-2mm}
\end{figure}

In Fig.\;\ref{fig10}, we show the first 100 singular values of the Loewner matrices computed by means of the three proposed partition techniques in (\ref{partition_types}) (left sub-figure). Additionally, the magnitude of the approximation errors on $[-1,1]$ for each of the four approximants is depicted (right sub-figure).

\begin{figure}[ht]		
	\begin{center}
		\includegraphics[scale=0.36]{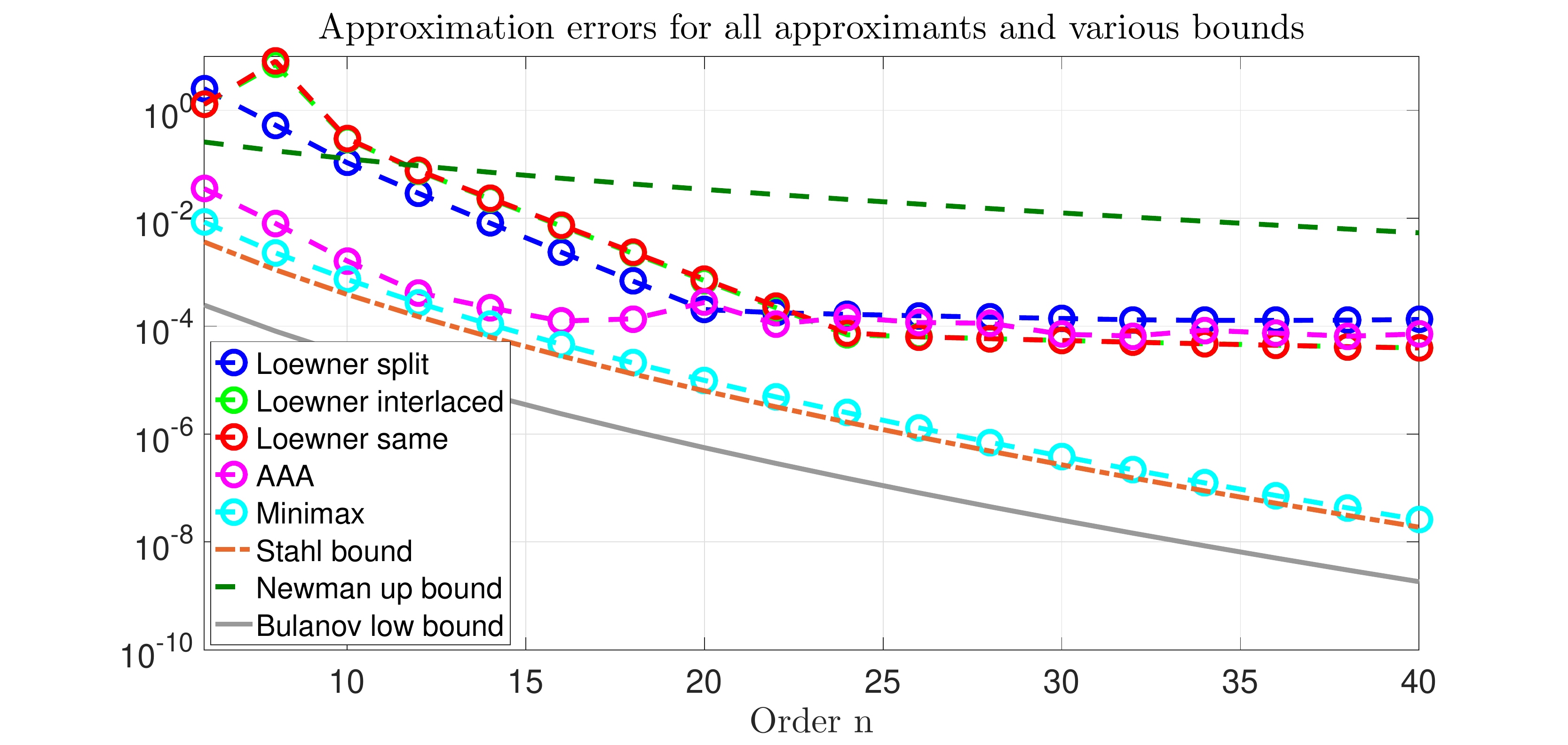}
	\end{center}
	\vspace{-3mm}
	\caption{The second case (logspace points). The maximum approximation errors for all methods and for different values of order $n$ + various classical bounds.}
	\label{fig11}
	\vspace{-2mm}
\end{figure}

Finally, in Fig.\;\ref{fig11} we depict the absolute value of the maximum approximation errors for all 4 methods and for various values of the rational functions' order (in the interval $[6,40]$). Additionally, some of the bounds mentioned in Section\;\ref{sec:approx} are also included. 

\subsection{Chebyshev points}

In this section we consider that the interpolation points are Chebyshev-type chosen within the interval $[-b,-a]\cup [a,b]$ for $a=2^{-10}$ and $b=1$. Let $n >0$ an integer and define positive Chebyshev points $p^{\rm Che}_1 < p^{\rm Che}_2 < \cdots < p^{\rm Che}_n$ as
\begin{equation}
p_k^{\rm Che} = \frac{1}{2}(a+b)+\frac{1}{2}(a-b)\cos\Big{(}\frac{(2k-1)\pi}{2n}\Big{)}, \ \ 1\leq k \leq n.
\end{equation}
Note that these points span a subinterval of the interval $(0,1)$. Additionally, for symmetry reasons, the negative points are also considered, i.e. $p_{-k}^{\rm Che} = - p_k^{\rm Che}$.

\begin{figure}[ht]		\includegraphics[scale=0.27]{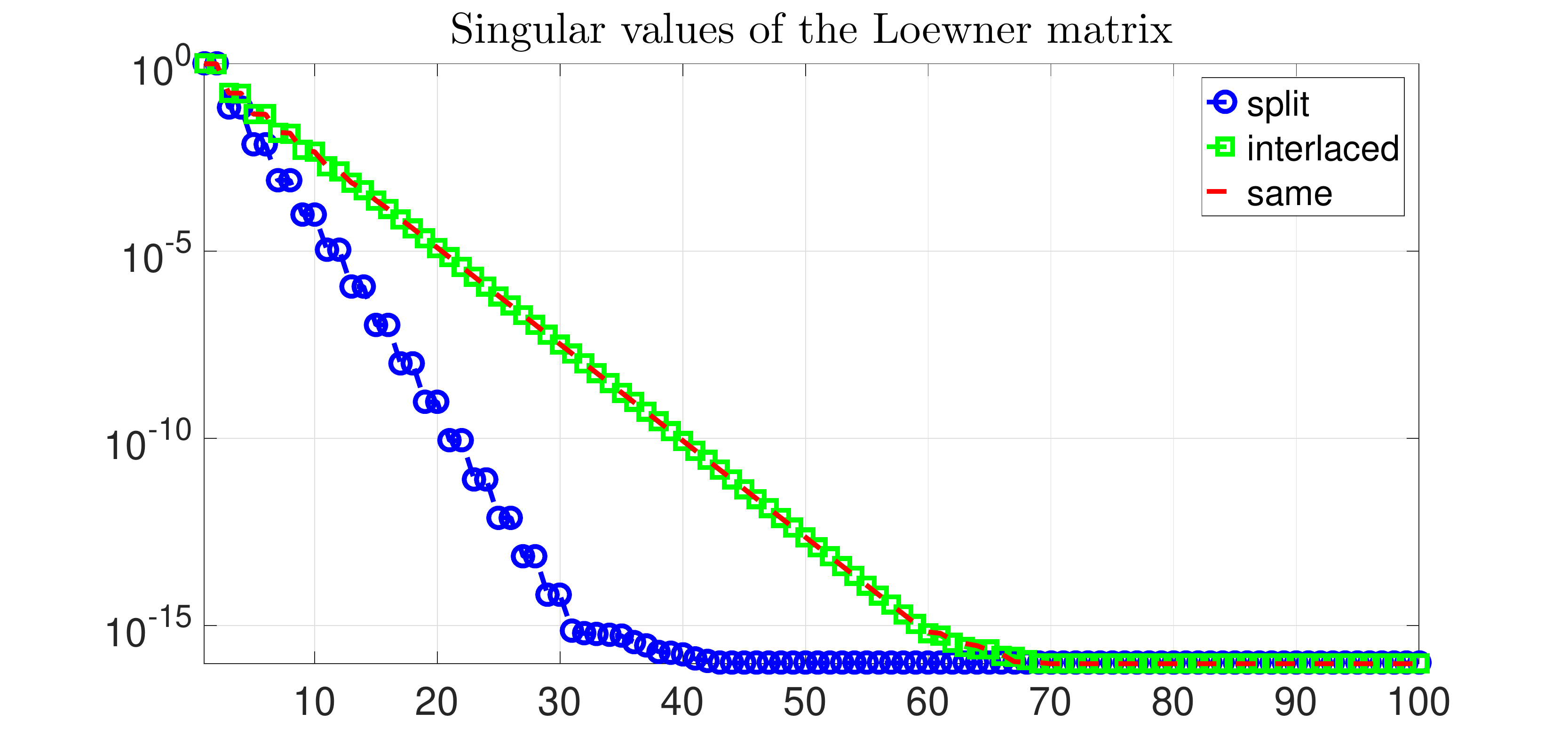}\hspace*{-8mm}	\includegraphics[scale=0.28]{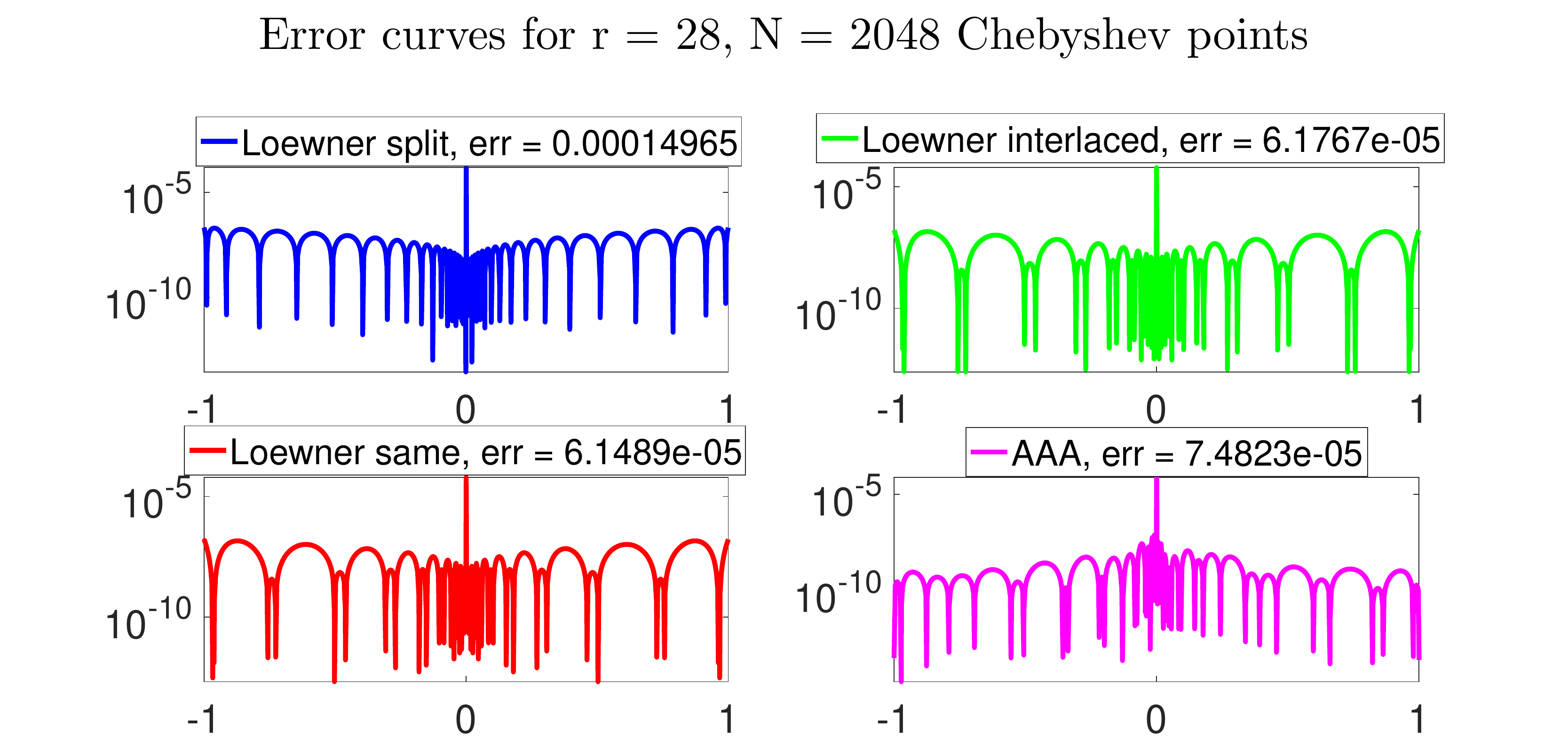}
	\vspace{-3mm}
	\caption{The third case (Chebyshev points).  Left: singular value decay of the Loewner matrices. Right: the magnitude of the error curves produced by order $r=28$ rational approximants.}
	\label{fig8}
	\vspace{-2mm}
\end{figure}

In Fig.\;\ref{fig8}, we show the singular value decay of the Loewner matrices computed by means of the three proposed partition techniques in (\ref{partition_types}) (left sub-figure). Additionally, the magnitude of the approximation errors on $[-1,1]$ for each of the four approximants is depicted (right sub-figure).

\begin{figure}[ht]		
	\begin{center}
		\includegraphics[scale=0.36]{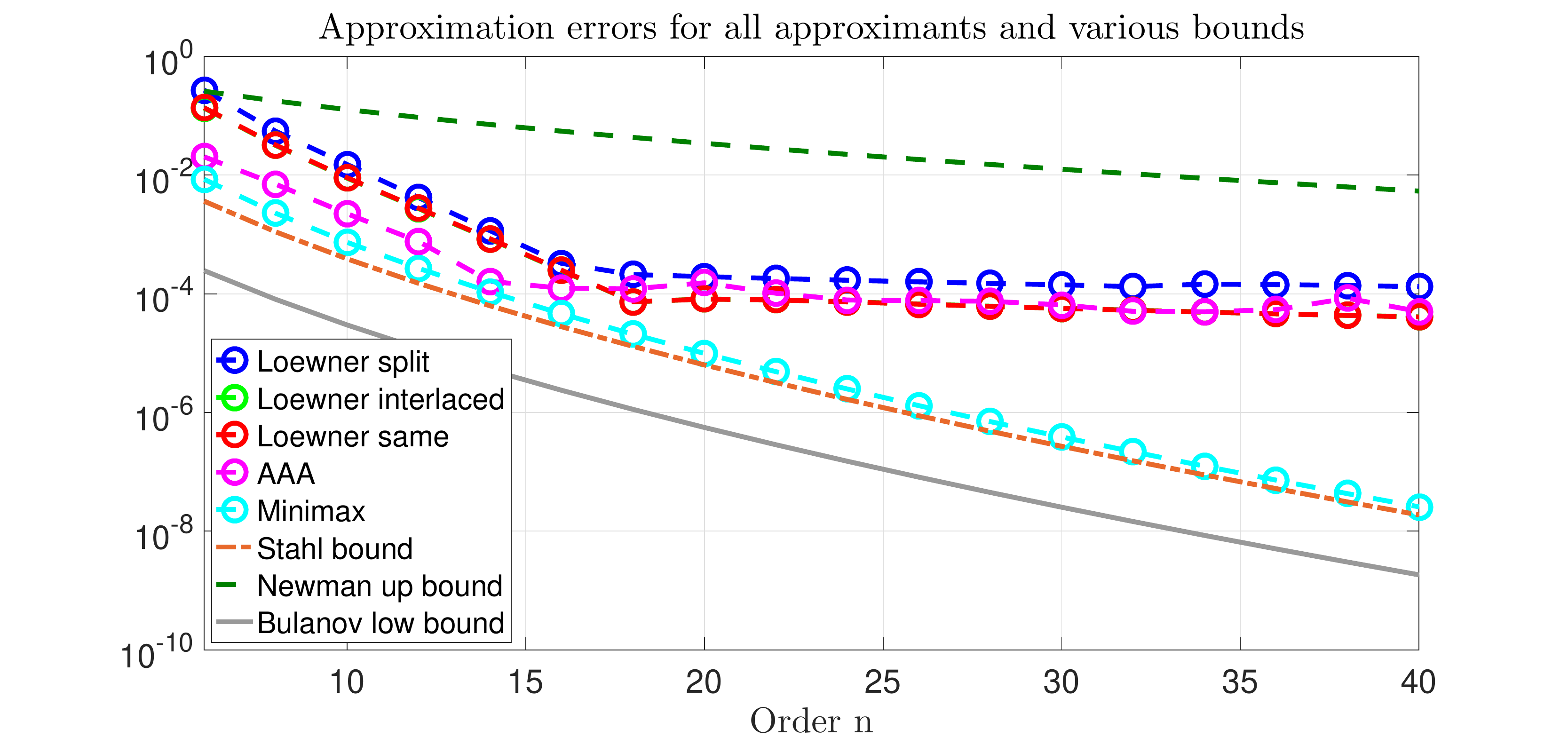}
	\end{center}
	\vspace{-3mm}
	\caption{The third case (Chebyshev points). The maximum approximation errors for all methods and for different values of order $n$ + various classical bounds.}
	\label{fig9}
	\vspace{-2mm}
\end{figure}

Finally, in Fig.\;\ref{fig9} we depict the absolute value of the maximum approximation errors for all 4 methods and for various values of the rational functions' order (in the interval $[6,40]$). Additionally, some of the bounds mentioned in Section\;\ref{sec:approx} are also included. 

\subsection{Zolotarev points}

In this section we consider that the given interpolation points are Zolotarev-type selected within the interval $[-b,-a]\cup [a,b]$. In 1877, Yegor Ivanovich Zolotarev explicitly found the best (minimax) rational approximation of the sign function 
\begin{equation}
{\rm sign}(x) = \begin{cases}
1, \ \ \text{if} \ x>0, \\
0, \ \ \text{if} \ x=0, \\
-1, \ \ \text{if} \ x<0. \\
\end{cases} 
\end{equation} 
For $0<a<b<1$, let $\ell = \frac{a}{b} <1$ and $\ell' = \sqrt{1- \ell^2} = \sqrt{1-\frac{a^2}{b^2}}$. Let $K'$ be the complete elliptic integral of the first kind for the complementary modulus and can be defined as
\begin{equation}
K' = \int_{0}^{\frac{\pi}{2}} \frac{d \theta}{\sqrt{1-(\ell')^2 \sin^2\theta}}. 
\end{equation}
Next, introduce $\rm{sn}(u;\ell')$ and $\rm{cn}(u;\ell')$ as the so-called Jacobi elliptic functions. The zeros and poles of the Zolotarev rational approximant of the sign function can be written in terms of these functions. For more details on such functions, see Chapter 5 in \cite{bookNA90}.

It can also be shown, i.e., as in \cite{nf16}, that the interpolation points of the Zolotarev approximant are explicitly given in terms of Jacobi elliptic functions. More precisely, these are the points in the interval $(a,b)$ where the sign function is exactly matched by the optimal Zolotarev approximant, i.e. where the value is exactly 1.

These points will be hence called "Zolotarev points" in what follows. They can be computed as follows; let $n >0$ to be a positive integer and $p_1^{\rm Zol} <p_2^{\rm Zol} < \cdots < p_n^{\rm Zol}$ be points such as
\begin{equation}
p_i^{\rm Zol} = \sqrt{a^2 \text{sn}^2(\frac{2i K'}{2n};\ell')+ b^2 \text{cn}^2(\frac{2i K'}{2n};\ell')}, \ \ \text{for} \ 1 \leq i \leq n.
\end{equation}
Note that these points span a subinterval of the interval $(0,1)$. Additionally, for symmetry reasons, the negative points are also considered, i.e. $p_{-i}^{\rm Zol} = - p_i^{\rm Zol}$.

\begin{figure}[ht]	
	\begin{center}
		\includegraphics[scale=1]{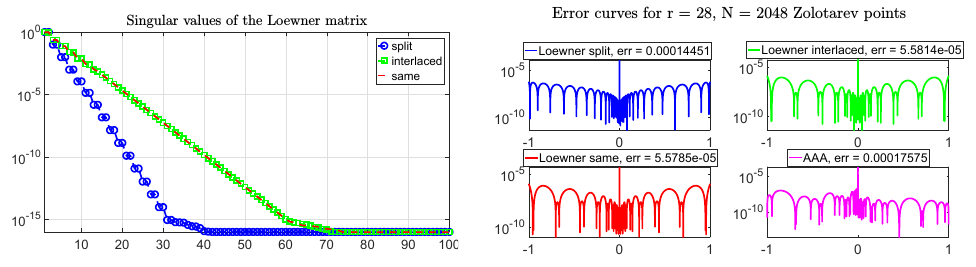}
		\end{center}
	\vspace{-3mm}
	\caption{The fourth case (Zolotarev points). Left: singular value decay of the Loewner matrices. Right: the magnitude of the error curves produced by order $r=28$ rational approximants.}
	\label{fig12}
	\vspace{-2mm}
\end{figure}

In Fig.\;\ref{fig12}, we show the first 100 singular values of the Loewner matrices computed by means of the three proposed partition techniques in (\ref{partition_types}) (left sub-figure). Additionally, the magnitude of the approximation errors on $[-1,1]$ for each of the four approximants is depicted (right sub-figure).

\begin{figure}[ht]		
	\begin{center}
		\includegraphics[scale=0.36]{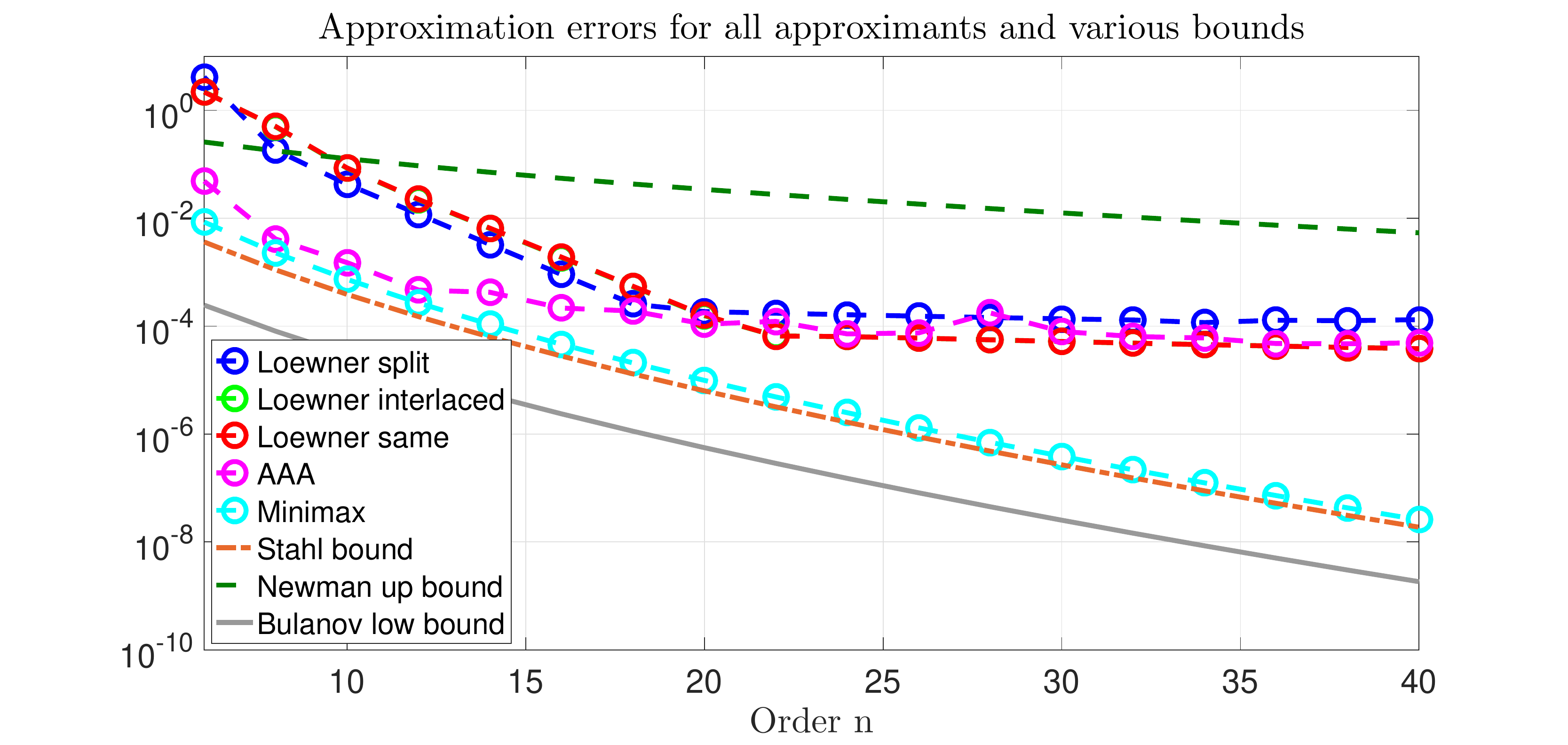}
	\end{center}
	\vspace{-3mm}
	\caption{The fourth case (Zolotarev points). The maximum approximation errors for all methods and for different values of order $n$ + various classical bounds.}
	\label{fig13}
	\vspace{-2mm}
\end{figure}

Finally, in Fig.\;\ref{fig13} we depict the absolute value of the maximum approximation errors for all 4 methods and for various values of the rational functions' order (in the interval $[6,40]$). Additionally, some of the bounds mentioned in Section\;\ref{sec:approx} are also included. 

In what follows, we also present some numerical results collected in two tables. This is done by comparing the four methods for the four types of sampling points (the case of "Newman points" in Section\;\ref{sec:Newman_pts} is deemed special and will be treated individually).

In Tab.\;\ref{table2} we show how the number of Loewner singular values that are bigger than the machine precision $2.2204 e-16$ vary with the type of sample points and partition schemes. Note that the fastest decay is encountered for "Loewner split" with linspace points, while the slowest decay is given by "Loewner same" with Zolotarev points (the number of singular values for the latter case is more than twice bigger than for the first case).

\begin{table}[h] 
	\begin{center}
		\begin{tabular}{ |p{2cm}||p{2.6cm}|p{3.4cm}|p{2.6cm}|  }
			\hline
			& Loewner split &  Loewner interlaced & Loewner same\\
			\hline
			Linspace   & 33    & 50 &   52\\
			Chebyshev & 37  & 65   & 64\\
			Logspace & 38 & 66 &  65\\
			Zolotarev & 38 & 68 &  68\\
			\hline
		\end{tabular}
	\end{center}
	\vspace{-4mm}
	\caption{Number of Loewner singular values bigger than machine precision.}
	\label{table2}
\end{table}

In Tab.\;\ref{table3} we collect the maximum approximation errors for all four methods and for the four types of sample points covered so far. The order of the rational approximants was kept constant, e.g., $r=28$. The highest error, i.e., $1.9920 e-04$ was obtained for "Loewner split" with linspace points while the lowest, i.e., $5.5785 e-05$ was obtained for "Loewner same" with Zolotarev points.

\begin{table}[h] 
	\begin{center}
		\begin{tabular}{ |p{2cm}||p{2.6cm}|p{3.4cm}|p{2.6cm}|p{2.6cm}|  }
			\hline
			& Loewner split &  Loewner interlaced& Loewner same & AAA\\
			\hline
			Linspace  & $1.9920 e-04$ & $9.8725 e-05$ & $7.9058 e-05$ & $1.0909 e-04$ \\
			Chebyshev & $1.4965 e-04$ & $6.1767 e-05$ & $6.1489 e-05$ & $7.4823 e-05$ \\
			Logspace & $1.9350 e-04$ & $1.9083 e-04$ & $1.9018 e-04$ & $1.5441 e-04$ \\
			Zolotarev & $1.4451 e-04$ & $5.5814 e-05$ & $5.5785 e-05$ & $1.7575 e-04$ \\
			\hline
		\end{tabular}
	\end{center}
	\vspace{-4mm}
	\caption{Approximation errors when the order was chosen $r=28$.}
	\label{table3}
\end{table}

\subsection{Newman points}\label{sec:Newman_pts}

Let $n >0$ an integer and $\alpha = e^{-\frac{\sqrt{n}}{n}} \in (0,1)$ and then define the Newton points $p^{\rm New}_1 < p^{\rm New}_2 < \cdots < p^{\rm New}_n$ as the powers of $\alpha$ from $n$ to $1$, i.e.
 \begin{equation}\label{points_Newman}
 p_1^{\rm New} = \alpha^n = e^{-\sqrt{n}}, \ \ p_2^{\rm New} = \alpha^{n-1} = e^{-\frac{(n-1)\sqrt{n}}{n}}, \ \ \ldots \ \ , p_n^{\rm New} = \alpha^1 = e^{-\frac{\sqrt{n}}{n}} .
 \end{equation}
These points span a subinterval of the interval $(0,1)$. Additionally, for symmetry reasons, the negative points are also considered, i.e. $p_{-k}^{\rm New} = - p_k^{\rm New}$. Note that the $2n$ points $p_k^{\rm New}$ for $-n \leq k \leq n$ are the interpolation points for constructing the order $(n,n)$ Newman approximant presented in Section\;\ref{sec:Newton} (the point $0$ was not included here).

For the next experiment, choose $N=256$ Newman points as described in \ref{points_Newman}.

\begin{figure}[ht]			\includegraphics[scale=1]{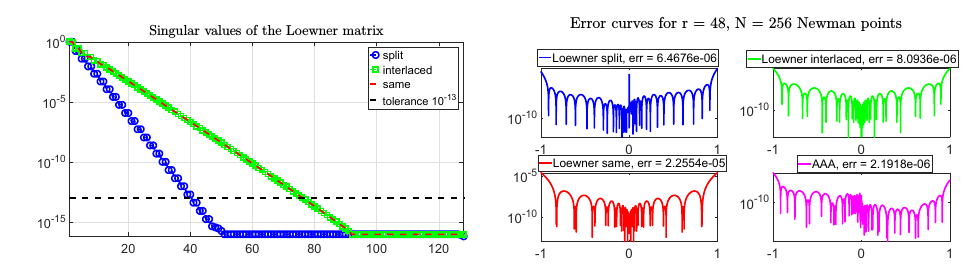}
	\vspace{-3mm}
	\caption{The fifth case (Newman points).  Left: singular value decay of the Loewner matrices. Right: the magnitude of the error curves produced by order $r=48$ rational approximants.}
	\label{fig14}
	\vspace{-2mm}
\end{figure}

 In Fig.\;\ref{fig14}, we show the 128 singular values of the Loewner matrices computed by means of the three proposed partition techniques in (\ref{partition_types}) (left sub-figure). Additionally, the magnitude of the approximation errors on $[-1,1]$ for each of the order  $r=48$ approximants is depicted (right sub-figure).

\begin{figure}[ht]		
	\begin{center}
		\includegraphics[scale=0.36]{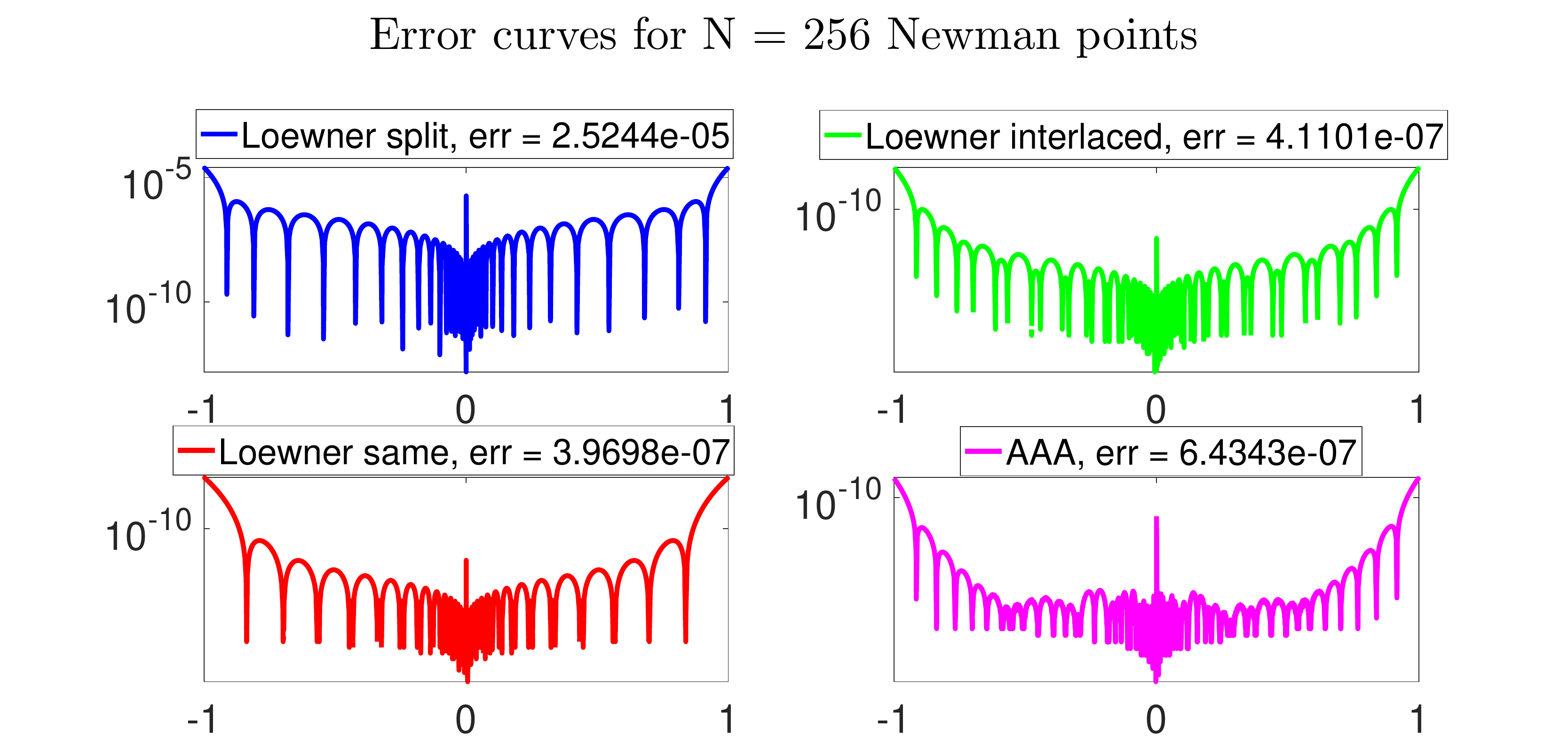}
	\end{center}
	\vspace{-3mm}
	\caption{The fifth case (Newman points).  The magnitude  of the error curves for the 4 rational approximants when $\delta = 1e-13$.}
	\label{fig15}
	\vspace{-2mm}
\end{figure}

Next, fix a tolerance of $\delta = 1e-13$ and truncate with respect to the number of singular values of the Loewner matrices that are larger than $\delta$. Hence, for the first partition technique, i.e. "Loewner split" we get $r_1 = 40$, while for the other two techniques we obtain $r_2 = r_3 = 76$. The reduced order for the $AAA$ method is chosen to be $r_4 = 76$.  Note that the result by Stahl presented in (\ref{Stahl_bound}) is $8 e^{-\pi \sqrt{76}} = 1.0203e-11$. In Fig.\;\ref{fig15}, we depict  the magnitude  of the approximation error curves for the all four rational approximants of orders $r_1$, $r_2$, $r_3$ and $r_4$.

\begin{figure}[ht]		
	\begin{center}
		\includegraphics[scale=0.36]{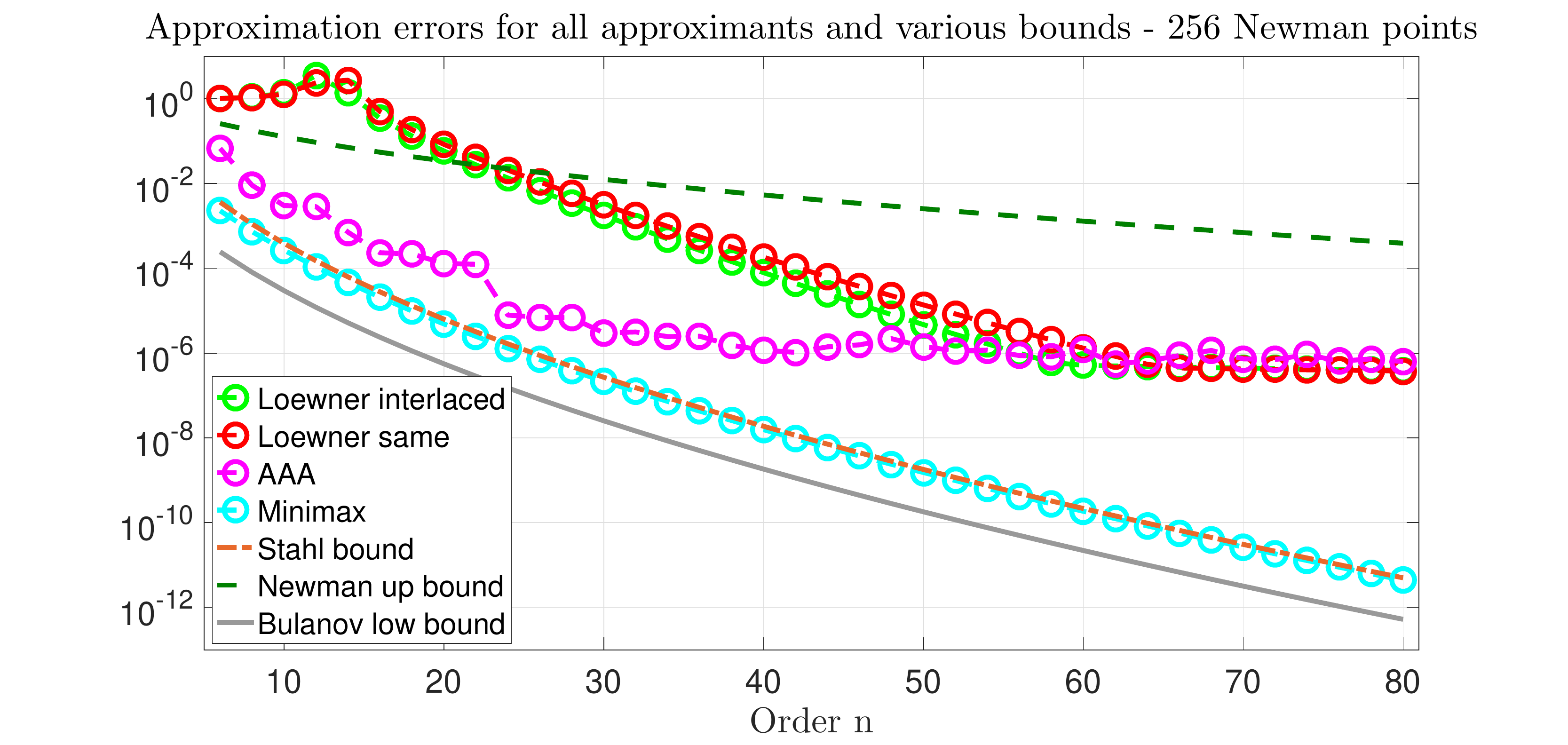}
	\end{center}
	\vspace{-3mm}
	\caption{The fifth case (Newman points). The maximum approximation errors for all methods and for different values of order $n$ + various classical bounds.}
	\label{fig16}
	\vspace{-2mm}
\end{figure}

Finally, in Fig.\;\ref{fig16} we depict the absolute value of the maximum approximation errors for all 4 methods and for various values of the rational functions' order (in the interval $[6,80]$). Additionally, some of the bounds mentioned in Section\;\ref{sec:approx} are also included. 

In Tab.\;\ref{table4} we show how the number of Loewner matrix singular values that are bigger than a tolerance value $\delta$ vary with the type of partition scheme. We observe that the decay reported for "Loewner split" is much faster than that of "Loewner interlaced" or "Loewner same".
 
\begin{table}[h] 
	\begin{center}
		\begin{tabular}{ |p{2.4cm}||p{2.6cm}|p{3.4cm}|p{2.6cm}|  }
			\hline
			& Loewner split &  Loewner interlaced & Loewner same \\
			\hline
			$\delta = 1e-9$  & $r_1 = 28$ & $r_2 = 54$ & $r_2 = 54$  \\
			$\delta = 1e-11$  & $r_1 = 34$ & $r_2 = 64$ & $r_3 = 64$ \\
			$\delta = 1e-13$ & $r_1 = 40$ & $r_2 = 76$ & $r_3 = 76$ \\
			$\delta = 1e-15$ & $r_1 = 46$ & $r_2 = 88$ & $r_3 = 88$ \\
			\hline
		\end{tabular}
	\end{center}
	\vspace{-4mm}
	\caption{Approximation orders and errors when the truncation tolerance value $\delta$ is varied.}
	\label{table4}
\end{table}

In Tab.\;\ref{table5} we collect the maximum approximation errors for all four methods and for different values of the tolerance $\delta$. The choice of $\delta$ determines the orders collected in Tab.\;\ref{table4}, for each of the four rational approximants.

Note that the highest errors were obtained for the "Loewner split" approximant. For the first 2 values of $\tau$, "Loewner interlaced" produced the lowest errors, while the last 2 values of $\tau$, "Loewner same" produced the lowest errors.

\begin{table}[h] 
	\begin{center}
		\begin{tabular}{ |p{2.4cm}||p{2.8cm}|p{3.4cm}|p{2.8cm}|p{2.8cm}|  }
			\hline
			& Loewner split &  Loewner interlaced& Loewner same & AAA\\
			\hline
			$\delta = 1e-9$  & $ 1.6347e-03$ & $1.6486e-06$ & $5.2023e-06$ & $ 1.1786e-06$ \\
			$\delta = 1e-11$  & $ 2.5377e-04$ & $4.8252e-07$ & $5.4574e-07$ & $ 6.5291e-07$ \\
			$\delta = 1e-13$ & $2.5244 e-05$ & $4.1101 e-07$ & $3.9698 e-07$ & $6.4343 e-07$ \\
			$\delta = 1e-15$ & $8.5655 e-06$ & $5.2722 e-07$ & $3.6259 e-07$ & $3.7328 e-07$  \\
			\hline
		\end{tabular}
	\end{center}
	\vspace{-4mm}
	\caption{Maximum approximation errors when the truncation tolerance value $\delta$ is varied.}
	\label{table5}
\end{table}

For the next experiment, choose eight times more points as for the previous experiment, e.g., $N=2048$ Newman points as given in \ref{points_Newman}. Note that $p_1^{\rm New} =  1.2664e-14$ and $p_{1024}^{\rm New} =  0.9692$. Hence, the left boundary of the (positive) interval in which the points are chosen is very close to 0. Additionally, the right boundary is close to 1, but not exactly 1.

\begin{figure}[ht]
	\begin{center}		\includegraphics[scale=0.36]{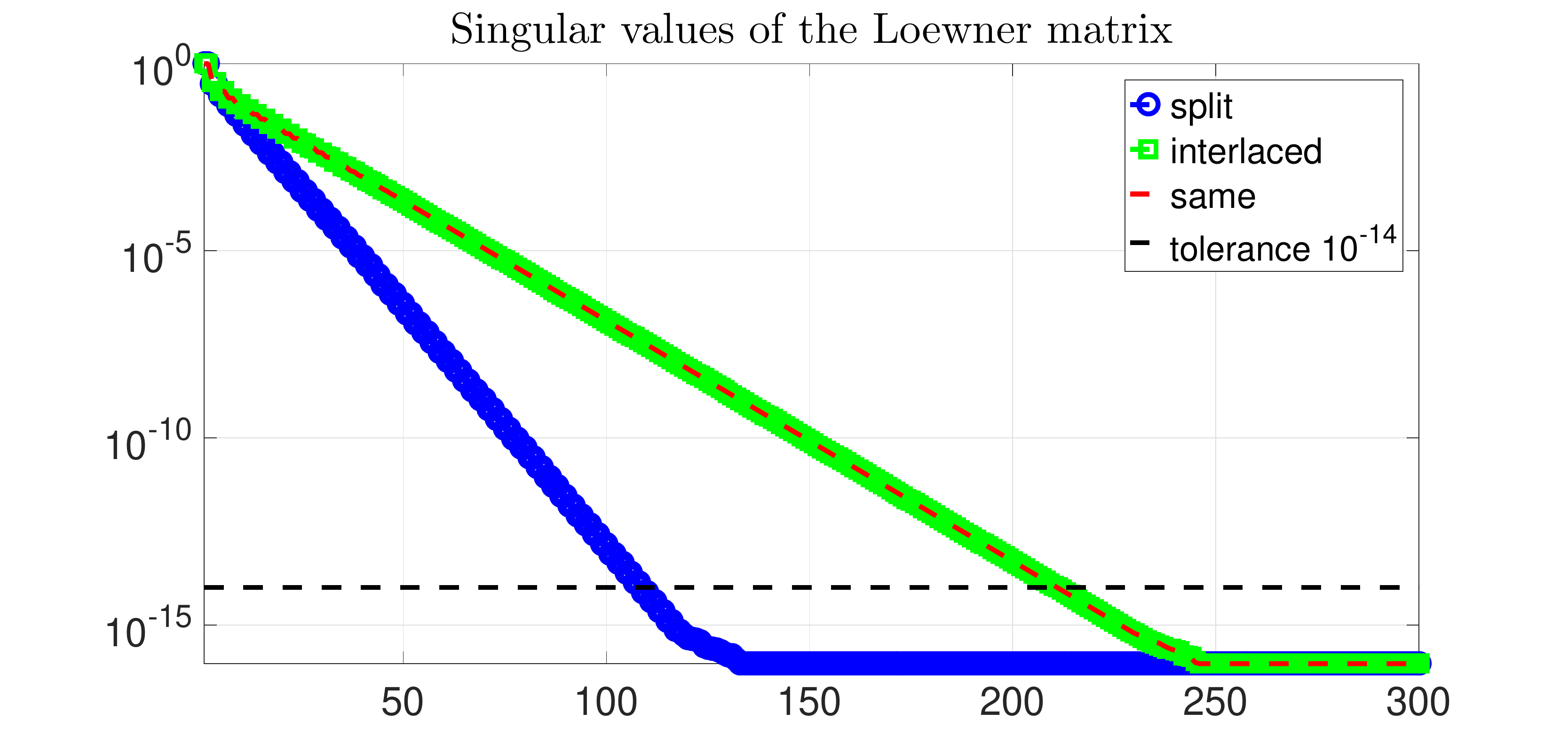}
	\end{center}
	\vspace{-3mm}
	\caption{The fifth case (2048 Newman points). Singular value decay of the Loewner matrices.}
	\label{fig17}
	\vspace{-2mm}
\end{figure}

In Fig.\;\ref{fig17}, we depict the first $300$ singular values of the Loewner matrices computed by means of the three proposed partition techniques in (\ref{partition_types}). Note that the dimension of these Loewner matrices is $1024 \times 1025$ (0 was again added in the right data set).

As observed in Fig.\;\ref{fig17}, the decay of singular values is quite slow (for all cases). This behavior is particularly noticed for "Loewner interlaced" and "Loewner same" partitioning cases. For each of these two variants, the machine precision is reached after around $250$ singular values. Choose a tolerance value $\tau = 1e-14$ and notice that a number of $210$ singular values are larger than this value (for the two cases).

\begin{figure}[ht]		
	\begin{center}
		\includegraphics[scale=0.36]{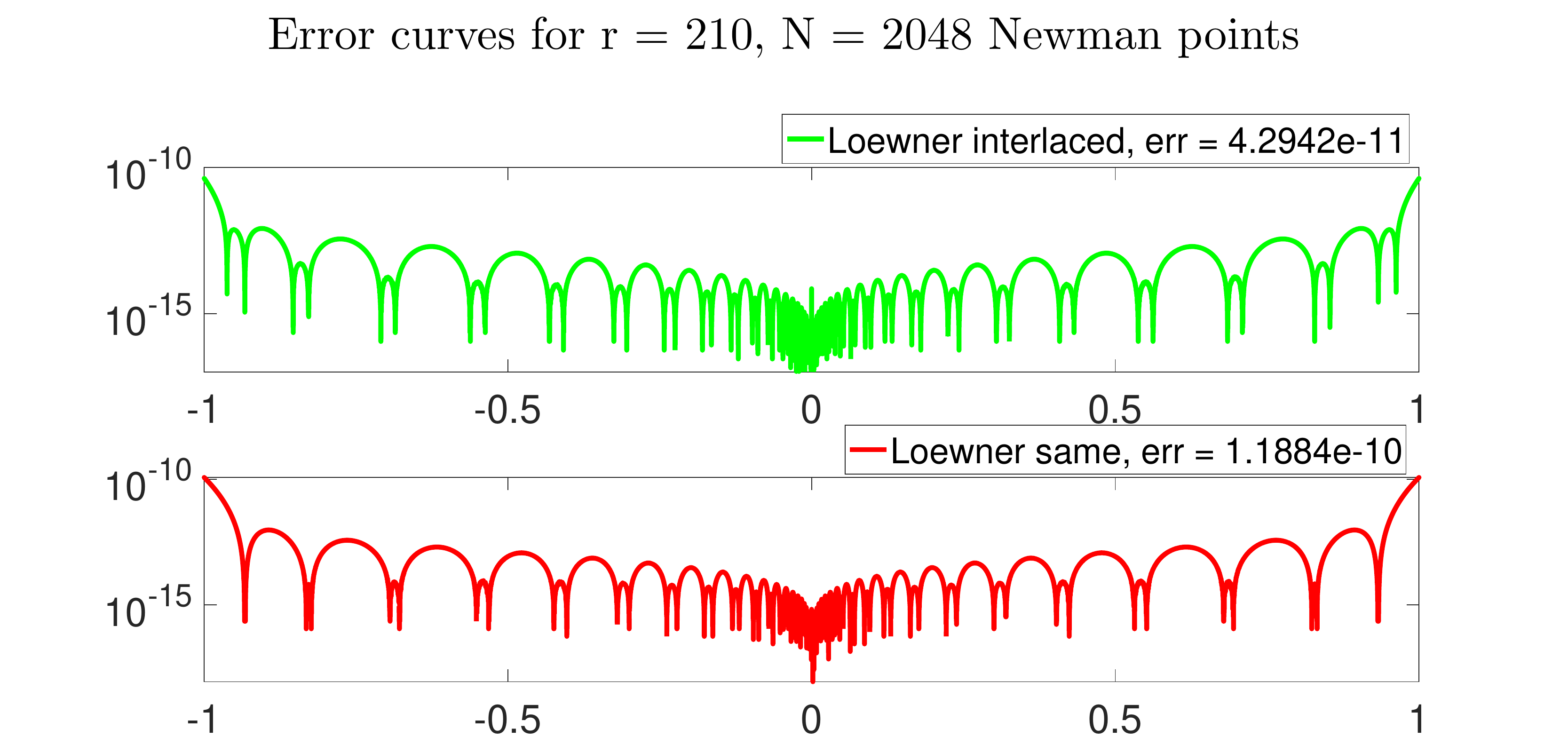}
	\end{center}
	\vspace{-3mm}
	\caption{The fifth case (2048 Newman points).  The magnitude  of the error curves for 2 rational approximants of order $r=210$.}
	\label{fig18}
	\vspace{-2mm}
\end{figure}

Next, based on the chosen tolerance value $\tau = 1e-14$, let $r=210$ be the corresponding reduction order for both of the reduced models computed with the Loewner procedure, i.e. "Loewner interlaced" and "Loewner same". In Fig.\;\ref{fig18}, we depict  the magnitude  of the approximation error curves for the two rational approximants of order $r$. Note that the maximum  approximation error in the case of "Loewner interlaced" drops below $1e-10$, i.e., reaching the value $4.2942 e-11$. It is to be denoted that, when $r=210$, it follows that the value from Strahl's bound in (\ref{Stahl_bound}) is $8 e^{-\pi \sqrt{210}} = 1.3533e-19$.

\section{An iterative Loewner framework procedure}\label{sec:iterate}

In this section we propose an iterative extension of the Loewner framework that can be used to increase the approximation quality of the Loewner rational approximant in (\ref{loew_fct}). The main idea is to add additional sample points and values to the original data set. The motivation for doing this is to decrease the overall maximum approximation error by improving the approximation at the newly added points. Through this procedure, the order of the Loewner rational approximant will remain constant (only the data set will increase in size).  

It is assumed that $(0,0)$ is not originally included in $\mathfrak{D}$. This is typically the pair that produces the maximum approximation error. Other pairs with high approximation errors are given by the extremes $(-1,1)$ and $(1,1)$. That is why, for the first two steps of the proposed procedure, precisely the pairs $\{(0,0), (-1,1), (1,1)\}$ will be added to $\mathfrak{D}$.

As before, we are given $N>0$ sampling points $\cT = \{\tau_1, \tau_2, \ldots,\tau_N\}$ together with $N$ function evaluations at these points denoted with $\{f_1, f_2, \ldots,f_N\}$. We have that $\tau_\ell \in [-b,-a] \cup [a,b]$ with $0 <a < b$ and $f_\ell = \vert \tau_\ell \vert,  \ \forall 1 \leq \ell \leq N$. The original set of measurement pairs is again denoted as in (\ref{data_Loew}) with $\mathfrak{D} = \{(\tau_\ell; \ \vert \tau_\ell \vert) \vert \ell =1,\cdots,N\}$.

In the newly proposed procedure, at each step, 2 more sampling points and values will be added to $\mathfrak{D}$. These extra points will be chosen similarly as for the AAA algorithm, i.e., by means of a Greedy procedure. More precisely, after step 2, the selected points are exactly the ones where the maximum deviation of the error function $\vert R_{\rm Loew}(x) - \vert x \vert \vert$ is attained on $(-1,0)$, and respectively on $(0,1)$. The pseudocode of the newly proposed algorithm is provided in Fig.\;(\ref{alg:LoewIter}).

\begin{figure}
	\noindent\fbox{%
		\parbox{.97\textwidth}{%
			\textbf{The Loewner iterative algorithm for approximating $\vert x \vert$}\\[10mm]
\noindent
\textbf{Input}: A set of measurements $\mathfrak{D} = \{(\tau_\ell; \ \vert \tau_\ell \vert) \vert \ell =1,\cdots,N\}$ corresponding to samples of the function $\vert x \vert$ on $[-b,-a] \cup [a,b]$, an integer $r > 0$, a positive tolerance value $\xi >0$. 

\noindent
\textbf{Output}: A rational approximant $\tilde{R}(x)$ of order $(r-1,r)$. \\[6mm]

\noindent
\textbf{Step 0}: Compute the Loewner rational approximant $R_{\rm Loew}(x)$ of order $r$ with respect to the initial data set $\mathfrak{D}$. Compute the values $\epsilon_R^{-}$ and $\epsilon_R^{+}$ and denote with $\rho \in (-1,1)$ the point where the maximum error occurs, i.e. $\vert R_{\rm Loew}(\rho) - \vert \rho \vert \vert = \max(\epsilon_R^{-},\epsilon_R^{+})$.\\

\noindent
\textbf{Step 1}: Include $(0,0)$ and $(\rho,\vert \rho \vert)$ to $\mathfrak{D}$. The new data set is $\mathfrak{D}^{(1)} = \mathfrak{D} \cup \{(0,0),(\rho,\vert \rho \vert)\} $. Compute the Loewner rational approximant $R_{\rm AAA}^{(1)}(x)$ of order $(r-1,r)$. \\

\noindent
\textbf{Step 2}: Include $(-1,1)$ and $(1,1)$ to $\mathfrak{D}$. The new data set is $\mathfrak{D}^{(2)} = \mathfrak{D}^{(1)} \cup \{(-1,1),(1,1)\} $. Compute the Loewner rational approximant $R_{\rm AAA}^{(2)}(x)$ of order $(r-1,r)$ and the overall maximum error on $[-1,1]$, i.e. $\epsilon_{R_{\rm Loew}^{(2)}}$. Let $m=2$.\\

\noindent
\textbf{While} $\Big{(}\epsilon_{R_{\rm Loew}^{(m)}} > \xi \Big{)}$ \\

\noindent
For $m \geq 2$, denote with $R_{\rm Loew}^{(m)}(x)$ the Loewner rational approximant of order $(r-1,r)$ at step $m$.  Compute the values $\epsilon_{R_{\rm Loew}^{(m)}}^{-}$ and $\epsilon_{R_{\rm Loew}^{(m)}}^{+}$ and denote with $\rho_m^{-} \in (-1,0), \rho_m^{+} \in (0,1)$ the points where these two maximum values are attained.\\

\noindent
The new data set is $\mathfrak{D}^{(m+1)} = \mathfrak{D}^{(m)} \cup \{(\rho_m^{-},\vert \rho_m^{-} \vert),(\rho_m^{+},\vert \rho_m^{+}\vert)\} $. Compute the Loewner rational approximant $R_{\rm Loew}^{(m+1)}(x)$ of order $(r-1,r)$ at step $m+1$. Compute the overall maximum error on $[-1,1]$, i.e. $\epsilon_{R_{\rm Loew}^{(m+1)}}$.\\

\noindent
If $ \epsilon_{R_{\rm Loew}^{(m+1)}} < \xi$, the algorithm stops. The sought after rational approximant of order $(r-1,r)$ is found, e.g., $\tilde{R}(x) = R_{\rm Loew}^{(m+1)}(x)$. If the error is still higher that $\xi$, continue the algorithm by adding two more measurement pairs to $\mathfrak{D}^{(m+1)}$. \\
 
\noindent
\textbf{End} }
}
\caption{The pseudocode of the Loewner iterative algorithm}
\label{alg:LoewIter}
\end{figure}

Next, we apply this algorithm with the following setup: $a = 2^{-10}, \ b = 1$, $N = 2048$ Chebyshev-type points, $r = 48$ and the tolerance value $\xi = 1e-07$. The partitioning scheme chosen is "same", both left and right data sets are identical. 

At step 0, we compute the maximum approximation error of the initial Loewner approximant: $\epsilon_{R_{\rm Loew}^{(0)}} =   1.1240e-04$. The goal is, after some iterations, to be able to decrease the maximum approximation error below a certain value, e.g. $\xi = 1e-07$.

\begin{figure}[ht]
	\begin{center}
	\includegraphics[scale=0.34]{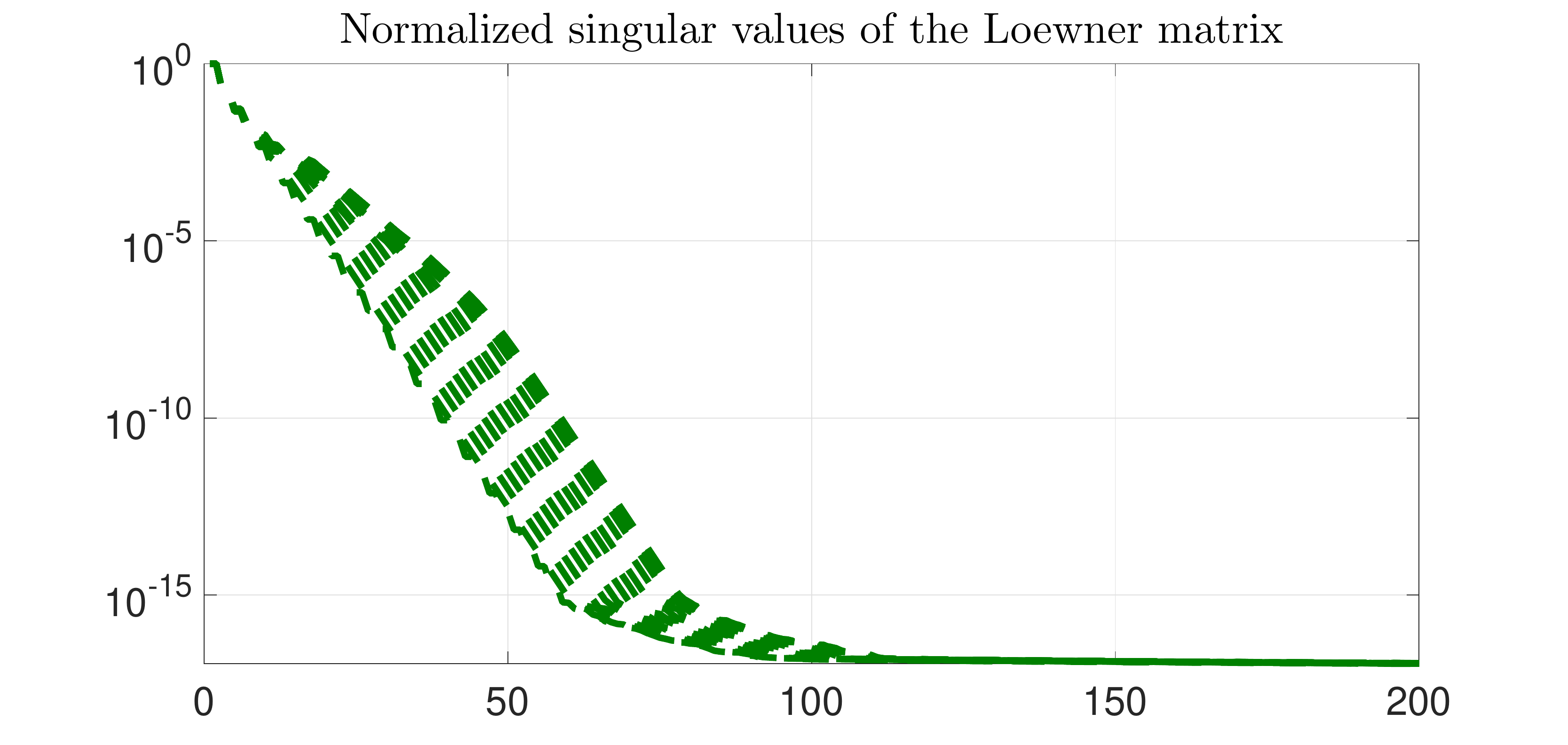}
	\end{center}
	\vspace{-5mm}
	\caption{The singular value decay of the Loewner matrices for all 22 steps.}
	\label{fig:plot19}
\end{figure}

We perform a total of $22$ steps. At each step, 2 sample points and values are added. The size of the Loewner matrix increases from $1024 \times 1024$ to $1048 \times 1048$. In Fig.\;\ref{fig:plot19} we depict the first 200 singular values of the 22 Loewner matrices, computed at each step.

\begin{figure}[ht]
	\begin{center}
	\includegraphics[scale=0.36]{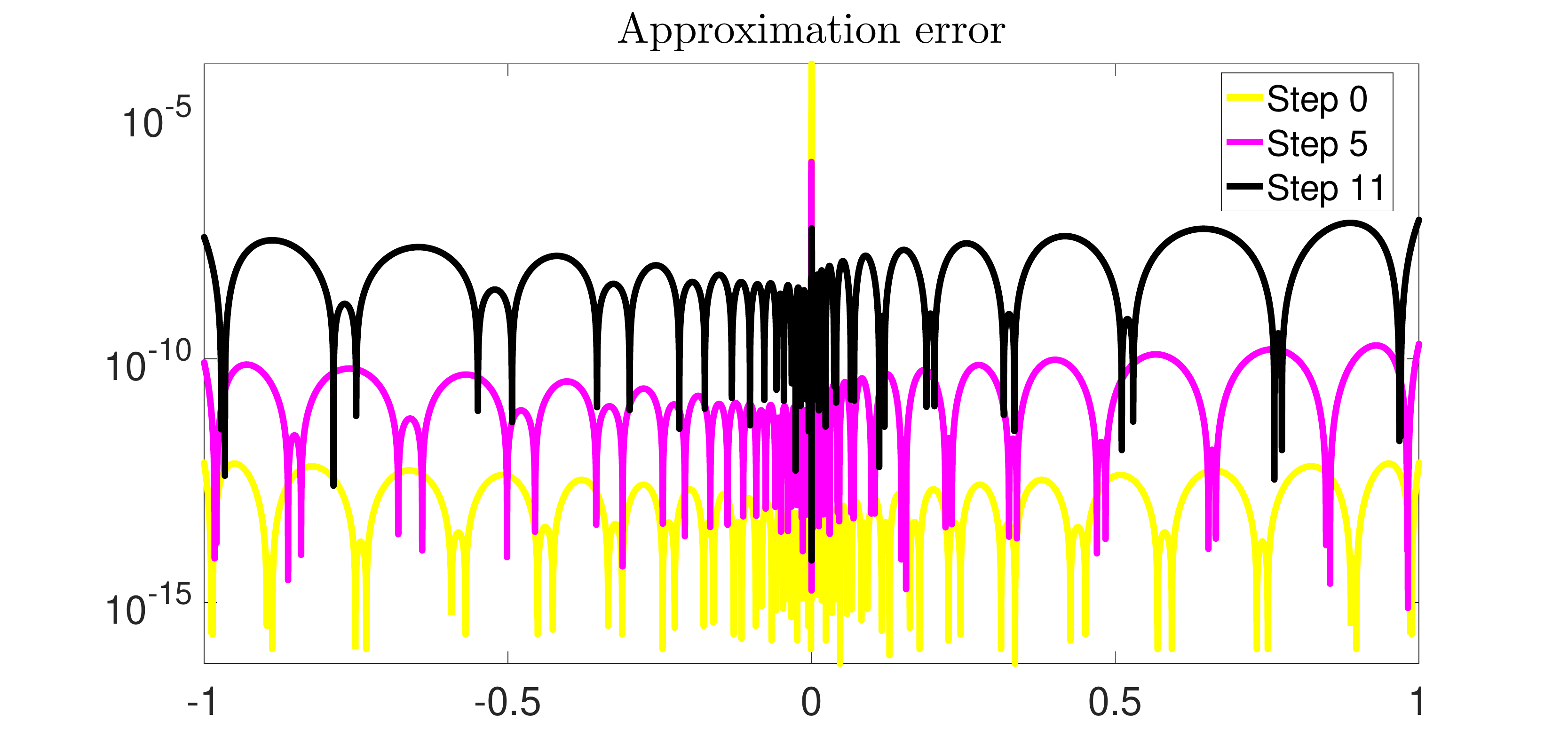}
	\end{center}
	\vspace{-5mm}
	\caption{The magnitude of the approximation error after 0, 5 and 11 steps.}
	\label{fig:plot20}
\end{figure}

\begin{figure}[ht]
	\begin{center}
	\includegraphics[scale=0.34]{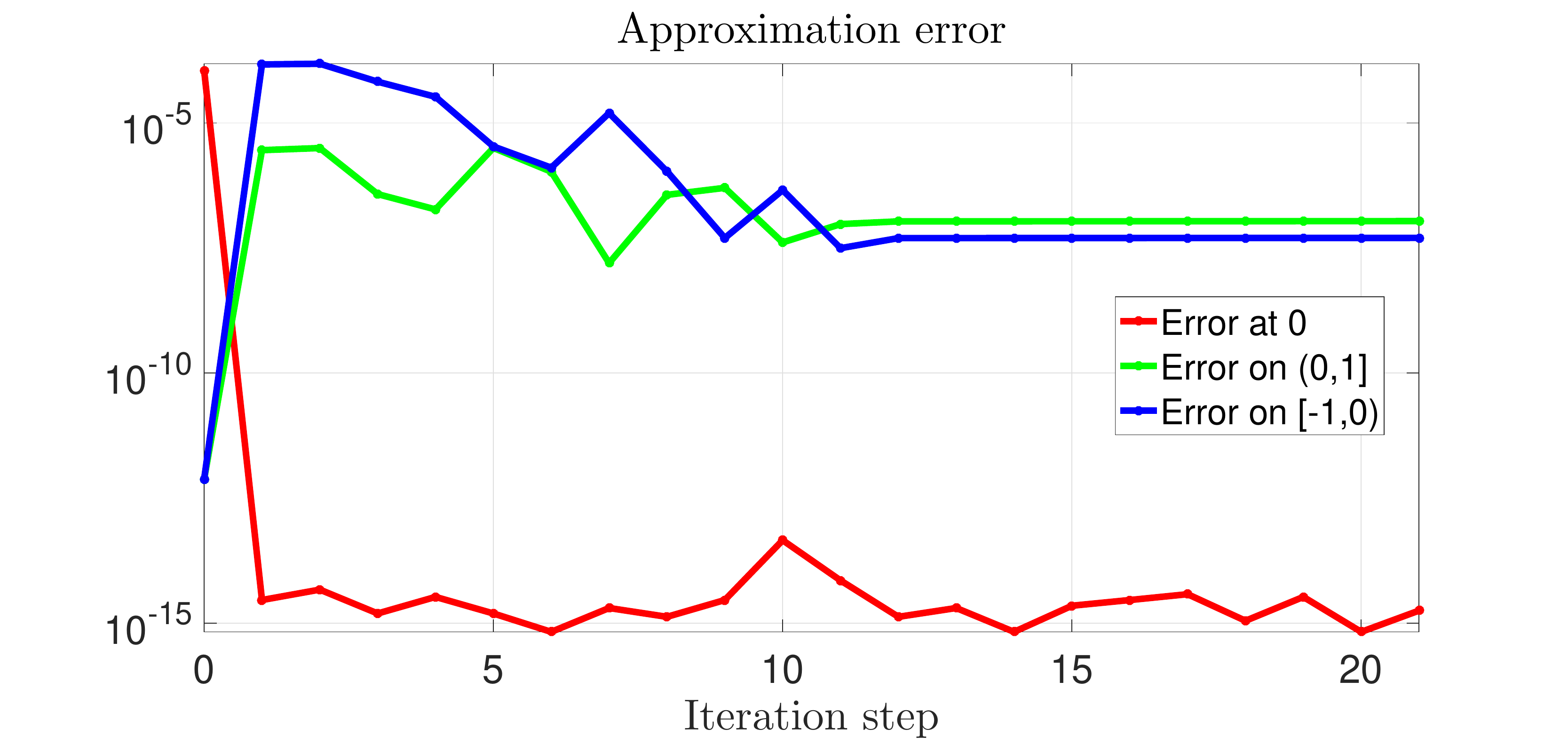}
    \end{center}
	\vspace{-5mm}
	\caption{The approximation error evaluated at 0 and on the intervals [0,1) $\&$ [-1,0).}
	\label{fig:plot21}
\end{figure}

\begin{figure}[ht]
    \begin{center}
	\includegraphics[scale=0.36]{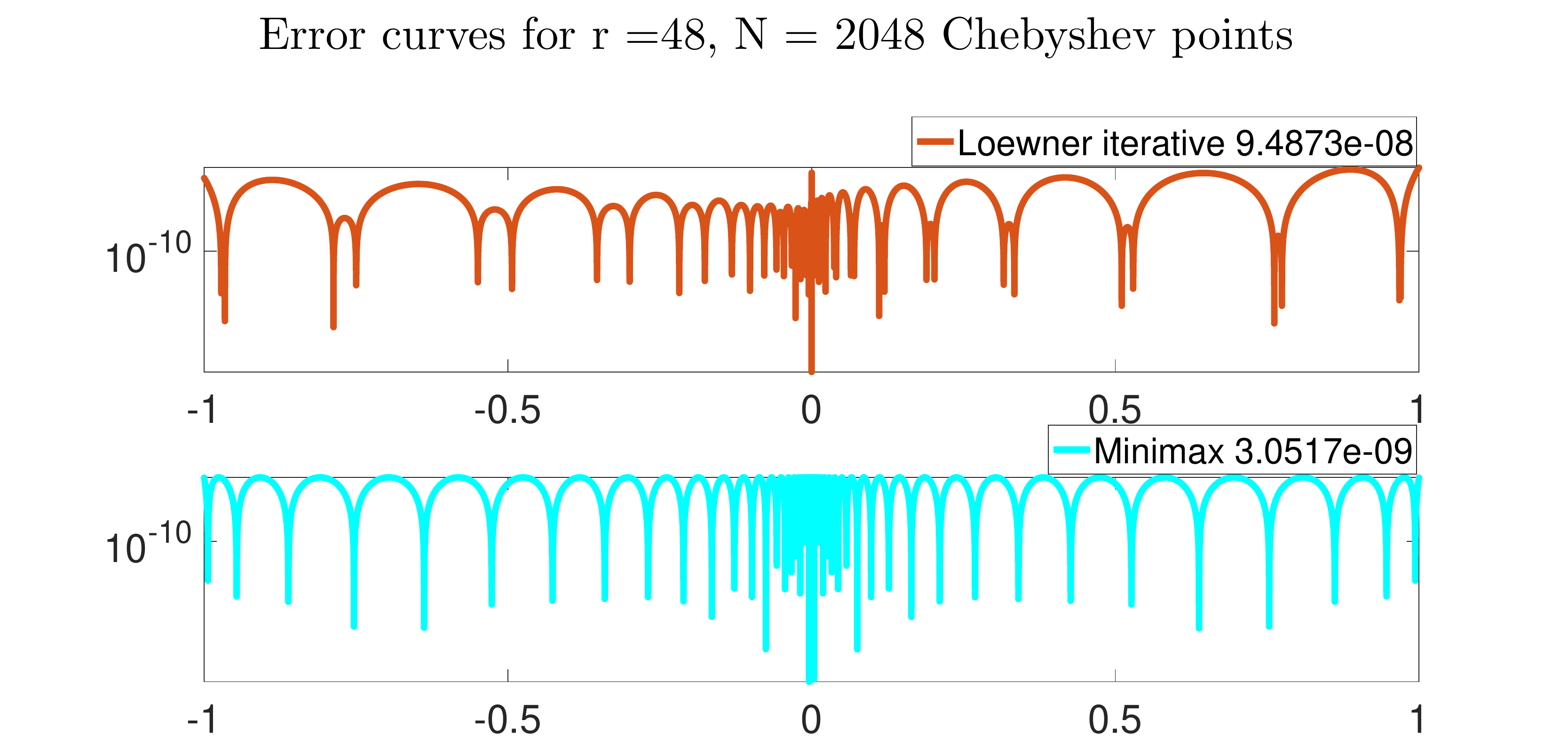}
	\end{center}
	\vspace{-5mm}
	\caption{The approximation error for Loewner iterative and minimax.}
	\label{fig:plot22}
\end{figure}

For each step $m$ with $0 \leq m \leq 21$, we compute the Loewner rational approximant $R_{\rm Loew}^{(m)}(x)$ as well as the maximum approximation error $\epsilon_{R_{\rm Loew}^{(m)}}$. In Fig.\;\ref{fig:plot20} we depict the magnitude of the approximation curves for three Loewner approximants $R_{\rm Loew}^{(m)}(x), \ m \in \{0,5,11\}$.
Next, for each step $m \in \{0,1,\ldots,21\}$ we record the approximation at the origin (0), as well as on the intervals $[-1,0)$ and $(0,1]$. The three curves are depicted in Fig.\;\ref{fig:plot21}. After step $m=11$, the curves are more or less flat.  We observe that, by imposing the tolerance $\xi = 1e-07$, the algorithm should have stopped at step $m=11$.
Finally, in Fig.\;\ref{fig:plot22}, we shop the magnitude of the approximation errors for the following two approximants:
\begin{enumerate}
	\item The first one, of order $(47,48)$, computed by means of the newly proposed Loewner iterative procedure: the maximum error is $ \epsilon_{R_{\rm Loew}^{(11)}} = 9.4873 e-08 < \xi$.
	\item The second one,  of order $(48,48)$, computed by means of the minimax algorithm in \cite{fntb18}: the maximum error is $ \epsilon_{R_{\rm Min}} = 3.0517 e-09$.
\end{enumerate} 
Note that, by computing the value of Stahl's result in (\ref{Stahl_bound}), we obtain $8e^{-\pi \sqrt{48}}= 2.8211 e-09$.

\begin{remark}
It is to be noted that the newly proposed algorithm in this section can easily be extended for approximating other functions, and it is by no means restricted for the example chosen in this study. The exact computation of the maximum error could be indeed challenging for non-rational functions. Nevertheless, for model order reduction purposes (approximating the rational transfer function of a large-scale system with a low order rational transfer function), this method can be used with no additional challenges. 
\end{remark}

\section{Conclusion}\label{sec:conc}

In this work we have proposed a detailed study that compares different data-driven methods  for approximating the absolute value function by means of rational functions. We also provided extensive numerical results that addressed issues such as using different type of measurements (sample points) or using different data splitting techniques. A detailed comparison to best rational approximants has also been included, together with carefully accounting for different results and bounds from classical approximation theory. Finally, we provided a novel iterative extension of the main method under investigation, i.e., the Loewner framework, in order to increase the overall approximation quality.

\newpage


\bibliographystyle{plain}

\bibliography{Gosea_Antoulas_Abs20}

\end{document}